\documentclass[11pt]{article}

\usepackage{amsmath, amssymb, amsfonts}
\usepackage{amsthm}
\usepackage{graphicx}
\usepackage{epstopdf}
\usepackage{algorithmic}
\usepackage{tabularx}
\usepackage{array}
\usepackage{booktabs}
\usepackage{mathtools}
\usepackage{pgfplots}
\pgfplotsset{compat=1.18}
\usepackage{tikz}
\usetikzlibrary{arrows.meta, positioning}
\usepackage{adjustbox}
\usepackage{multirow}
\usepackage{subcaption}
\usepackage{xspace}
\usepackage[colorlinks=true,linkcolor=blue,citecolor=blue,urlcolor=blue]{hyperref}
\usepackage{cleveref}

\usepackage[a4paper,top=1in,bottom=1in,left=0.9in,right=0.9in]{geometry}

\newcommand{\yauyau}{Yau--Yau\xspace}

\newtheorem{theorem}{Theorem}[section]
\newtheorem{lemma}[theorem]{Lemma}
\newtheorem{definition}[theorem]{Definition}

\newtheorem{assumption}[theorem]{Assumption}
\newtheorem{remark}[theorem]{Remark}

\numberwithin{equation}{section}

\title{An Improved Yau--Yau Algorithm for High Dimensional Nonlinear Filtering Problems}

\author{%
  Shing-Tung Yau\thanks{Yau Mathematical Sciences Center, Tsinghua University, Beijing 100084, China; 
  Beijing Institute of Mathematical Sciences and Applications (BIMSA), Beijing 101408, China. 
  Email: \texttt{styau@tsinghua.edu.cn}}%
  \and
  Yi-Shuai Niu\thanks{Beijing Institute of Mathematical Sciences and Applications (BIMSA), Beijing 101408, China. 
  Email: \texttt{niuyishuai@bimsa.cn}}
}

\date{}

\begin{document}

\maketitle

\begin{abstract}
Nonlinear state estimation under noisy observations is rapidly intractable as system dimension increases. We introduce an improved Yau–Yau filtering framework that breaks the curse of dimensionality and extends real-time nonlinear filtering to systems with up to thousands of state dimensions, achieving high-accuracy estimates in just a few seconds with rigorous theoretical error guarantees. This new approach integrates quasi‐Monte Carlo low‐discrepancy sampling, a novel offline-online update, high‐order multi‐scale kernel approximations, fully log‐domain likelihood computation, and a local resampling-restart mechanism, all realized with CPU/GPU‐parallel computation. Theoretical analysis guarantees local truncation error \(O(\Delta t^2 + D^*(n))\) and global error \(O(\Delta t + D^*(n)/\Delta t)\), where \(\Delta t\) is the time step and \(D^*(n)\) the star-discrepancy. Numerical experiments—spanning large‐scale nonlinear cubic sensors up to 1000 dimensions, highly nonlinear small‐scale problems, and linear Gaussian benchmarks—demonstrate sub‐quadratic runtime scaling,  
sub‐linear error growth, and excellent performance that surpasses the extended and unscented Kalman filters (EKF, UKF) and the particle filter (PF) under strong nonlinearity, while matching or exceeding the optimal Kalman–Bucy filter in linear regimes. By breaking the curse of dimensionality, our method enables accurate, real‐time, high‐dimensional nonlinear filtering, opening broad opportunities for applications in science and engineering.
\end{abstract}




\section{Introduction}
The problem of nonlinear filtering, estimating the state-space variable \(x_t\in\mathbb{R}^r\) of a continuous-time stochastic system from noisy observations, arises in diverse applications, including spacecraft navigation \cite{CoxBrading2000,Schmidt1981}, quantitative finance \cite{Wells1996}, target tracking \cite{Aidala1979}, weather forecasting \cite{Cassola2012}, and autonomous control \cite{Thrun2002}, among others. In the linear–Gaussian setting, the Kalman filter was first formulated for discrete time by Kalman \cite{Kalman1960}, and subsequently extended to continuous time by Kalman and Bucy \cite{KaBu}, which provides an exact finite‐dimensional solution.
Its extensions such as the Extended Kalman Filter \cite{Jazwinski1970} (EKF) and the Unscented Kalman Filter \cite{JulierUhlmann1997} (UKF) remain effective under mild nonlinearity \cite{FlMi}.  For the fully nonlinear setting, early theoretical foundations include Kushner’s dynamical equations for nonlinear filtering \cite{KuDE}, Kushner’s discrete-state approximations \cite{Ku} and Zakai’s unnormalized Stochastic Partial Differential Equation (SPDE) for the conditional density \cite{Za}, later extended by Duncan \cite{Du} and Mortensen \cite{Mo} to the continuous-time Duncan-Mortensen-Zakai (DMZ) equations.  Pardoux and collaborators then established existence, uniqueness, and regularity for these infinite-dimensional DMZ equations \cite{Pa1,Pa2,Pa3}, but their infinite-dimensional nature precludes direct solution.  

Practical algorithms approximate the infinite-dimensional DMZ problem via tractable finite representations. Spectral methods expand the density in orthogonal bases \cite{LMR}, controlling approximation error at the cost of exponentially many modes in high dimensions. Sequential Monte Carlo (particle) filters approximate the DMZ equation with weighted samples \cite{DJP}, but suffer from weight degeneracy and exponential computational cost as the state dimension \(r\) grows. Quasi–Monte Carlo (QMC) sampling with low-discrepancy sequences (e.g., Halton \cite{Halton1960} and Sobol \cite{Sobol1967}) offers substantial variance reduction, partially mitigating this “curse of dimensionality.” Operator-splitting methods have also been applied to the Zakai SPDE \cite{HoWo,BGR1,BGR2,GyKr,FlLe}, but they introduce an irreducible splitting error due to non-commutativity of the decomposed operators, which can be amplified in stiff or high-dimensional settings.

Independently, Yau \& Yau (2000, 2008) proposed a conceptually different memory-less, real-time filtering strategy, now known as the \yauyau filter \cite{Ya-Ya-I,Ya-Ya-II}. This method exploits an offline-online decomposition: the offline stage computes the deterministic (noise-free) propagator in advance, while the online stage performs simple multiplicative updates using incoming observations. Subsequent implementations by Yueh, Lin \& Yau \cite{YuehLinYau2014} validated the \yauyau filter in moderate dimensions ($r=2$). More recently, Chen et al. \cite{ChenSunTaoYau2025} introduced a deep learning-based \yauyau filter, in which an offline‐trained recurrent neural network (RNN) replaces the propagation‐and‐update steps to achieve real‐time filtering for \(r\le100\) with polynomial cost growth.  However, the extensive offline phase, which involves high-dimensional trajectory generation and training the RNN over thousands of epochs, takes hours of GPU computation, and no experiments beyond \(r=100\) have been reported, limiting applicability in higher dimensions. So far, scaling beyond a few hundred dimensions for nonlinear filtering has remained intractable.

In this work, we present an \textbf{improved Yau–Yau filtering framework} that breaks the curse of dimensionality and extends real-time nonlinear filtering to systems with up to thousands of state dimensions, achieving high-accuracy estimates in just a few seconds with rigorous error and stability guarantees. This breakthrough is built upon the integration of the following key innovations: (i) \textbf{Low‐discrepancy QMC sampling} for high‐dimensional state coverage; (ii) Novel \textbf{offline-online} update with provable local truncation error \(O(\Delta t^2)\) and global error bound \(O(\Delta t)\), where \(\Delta t\) is the time step; (iii) \textbf{High‐order, multi‐scale kernel approximations} of the Kolmogorov propagator; (iv) \textbf{Fully log‐domain likelihood updates} ensuring numerical stability; (v) \textbf{Local resampling-restart mechanism} to prevent weight collapse;
(vi) \textbf{CPU/GPU‐parallel implementations} to accelerate both offline and online stages in high dimensions.

Our theoretical analysis shows local and global error bounds of \(O(\Delta t^2 + D^*(n))\) and \(O(\Delta t + D^*(n)/\Delta t)\), respectively, where \(D^*(n)\) is the star-discrepancy of the QMC sequence.  Numerical experiments—including large‐scale tests up to $r=1000$ state dimension, highly nonlinear small‐scale systems, and linear‐Gaussian benchmarks—demonstrate sub‐quadratic runtime scaling (\(\sim r^{1.2}\)), sub‐linear error growth, and accuracy superior to EKF/UKF and PF in nonlinear regimes, while matching or exceeding the Kalman–Bucy in linear regimes.  These results effectively break the curse of dimensionality, making real-time nonlinear filtering practical to previously intractable, high‐dimensional applications.

\section{Improved Yau–Yau Algorithm}\label{sec:improved-Yau–Yau-algorithm}
\subsection{Continuous-Time Nonlinear Filtering Problem}
Consider the continuous-time nonlinear filtering problem
\begin{equation}
\label{eq:NonlinearFilter}
\begin{cases}
 \mathrm{d}X_t = f(X_t)\,\mathrm{d}t+ \mathrm{d}V_t, & X_0 = x_0,\\[1mm]
 \mathrm{d}Y_t = h(X_t)\,\mathrm{d}t+ \mathrm{d}W_t, & Y_0 = 0,\quad t\in[0,T],
\end{cases}
\end{equation}
where:
\begin{itemize}
  \item $X_t\in\mathbb{R}^r$ is the $r$-dimensional hidden state with initial value $x_0$.
  \item $f:\mathbb{R}^r\to\mathbb{R}^r$ is the drift (system dynamics) function.
  \item $V_t$ is an $r$-dimensional standard Brownian motion representing process noise.
  \item $Y_t\in\mathbb{R}^m$ is the $m$-dimensional observation with $Y_0=0$.
  \item $h:\mathbb{R}^r\to\mathbb{R}^m$ is the observation function.
  \item $W_t$ is an $m$-dimensional standard Brownian motion representing observation noise.
  \item $V_t$ and $W_t$ are independent.
\end{itemize}
All stochastic integrals are interpreted in the Itô sense.

Let $\sigma(x,t)$ be the unnormalized conditional probability density of $X_t$ given the observation history $\mathcal{Y}_t=\{y_s:0\le s\le t\}$.  It is well known (see \cite{Du,Mo,Za}) that \(\sigma(x,t)\) satisfies the Duncan–Mortensen–Zakai (DMZ) equation:
\begin{equation}
    \label{eq:DMZ}
\mathrm{d}\sigma(x,t)
= \mathcal{L}^*\,\sigma(x,t)\,\mathrm{d}t
+ \sigma(x,t)\,h(x)^{\top}\mathrm{d}Y_t,
\quad
\sigma(x,0)=\sigma_0(x),
\end{equation}
with the adjoint operator
\begin{equation}
    \label{eq:L*}
\mathcal{L}^* \sigma
= \frac{1}{2}\,\Delta\sigma
- \nabla\cdot \bigl(f(x)\,\sigma\bigr).
\end{equation}
Here,
\[
\Delta\sigma
= \sum_{i=1}^r \frac{\partial^2 \sigma}{\partial x_i^2}
\]
is the Laplace operator modeling diffusion, and
\[
\nabla\cdot \bigl(f(x)\,\sigma\bigr)
= \sum_{i=1}^r \frac{\partial}{\partial x_i}\!\bigl(f_i(x)\,\sigma(x,t)\bigr)
\]
is the divergence of the vector field \(f(x)\,\sigma\), representing the effect of the drift \(f\) on the density. The normalized conditional density is then
\[\frac{\sigma(x,t)}{\int \sigma(x,t)\mathrm{d}x}.\]

\subsection{Improved Yau–Yau Algorithm}
Consider solving the DMZ equation 
\[
\mathrm{d}\sigma(x,t)
= \mathcal{L}^*\sigma(x,t)\,\mathrm{d}t
+ \sigma(x,t)\,h(x)^{\top}\,\mathrm{d}Y_t,
\quad \sigma(x,0)=\sigma_0(x),
\]
where
\[
\mathcal{L}^*\sigma
= \frac{1}{2}\,\Delta\sigma
- \nabla\cdot \bigl(f(x)\,\sigma\bigr).
\]
In general, the DMZ equation \eqref{eq:DMZ} admits no closed‐form solution. Yau and Yau \cite{Ya-Ya-II} made a breakthrough in this problem by introducing a two-stage (online-offline stages) numerical algorithm, namely \yauyau algorithm, in which the
computationally expensive procedure of numerically solving the Kolmogorov forward equation
\[\frac{\partial \sigma}{\partial t} (x,t) = \mathcal{L}^*\sigma(x,t) - \frac{1}{2}\|h(x)\|^2 \sigma(x,t)\]
can be done offline for all \(t\) on a predefined mesh; and during the online stage, the solution is updated at observation time $\tau_k$ by the exponential transform
    \[
\tilde\sigma(x,\tau_k)
= \exp\!\Bigl\{h(x)^{\top}\Delta y_k\Bigr\}\,\sigma(x,\tau_k),  \quad \Delta y_k = y_{\tau_k}-y_{\tau_{k-1}}.
\]

In this section, we introduce an improved \yauyau framework based on low-discrepancy sequences for state-space sampling and a novel offline-online \yauyau scheme.

\subsubsection{Novel Offline-Online \yauyau Scheme with Low‐Discrepancy Sampling}

We begin by discretizing both time and state space:

\begin{itemize}
  \item \textbf{Time discretization:} Partition the interval \([0,T]\) into \(K\) equal subintervals
    \[
      0 = \tau_0 < \tau_1 < \cdots < \tau_K = T,
      \quad 
      \Delta t = \frac{T}{K}.
    \]
  \item \textbf{State sampling:} Generate \(n\) quasi‐uniform points 
    \(\{x_i\}_{i=1}^n \subset [-R,R]^r\) 
    via quasi‐Monte Carlo (QMC) low‐discrepancy sequences (e.g.\ Halton, Sobol, or Faure) or Latin hypercube sampling (LHS); 
    see Appendix~\ref{sec:low-discrepancy-sequence-sampling}.
\end{itemize}

To approximate the solution of DMZ \eqref{eq:DMZ} over each time interval \([\tau_{k-1},\tau_k]\) with step size $\Delta t$, we split the DMZ generator into two operators
\[
  A \sigma = \mathcal{L}^*\,\sigma,
  \qquad
  B \sigma = \sigma\,h(x)^\top \,\dot Y(t),
\]
and propose the following offline-online update:
\[
  \sigma_k
  = M(\Delta y_k)\;\exp\!\bigl(A\,\Delta t\bigr)\,\sigma_{k-1},
\]
where $\sigma_k$ is the discrete approximation of $\sigma(x,\tau_{k})$ and 
\[
  \Delta y_k = y_{\tau_k}-y_{\tau_{k-1}},
  \quad
  M(\Delta y_k)
  = \exp\left\{ h(x)^{\top}\Delta y_k
    - \tfrac12\,\|h(x)\|^2\,\Delta t\right\}.
\]
This scheme leverages the precomputation of \(\exp(A\,\Delta t)\) offline and applies the simple multiplicative update \(M(\Delta y_k)\) online, while the use of low‐discrepancy points \(\{x_i\}\) ensures an efficient computation of $\exp(A\Delta t)$ describing below.

\paragraph{Offline stage (Prediction)}
    The precomputation of \(\exp(A\,\Delta t)\) offline is equivalent to solve the pure‐dynamics Kolmogorov forward equation
    \[\begin{cases}
    \displaystyle \frac{\partial \rho(x,t)}{\partial t}
      = \mathcal{L}^*\,\rho(x,t),
      & t\in[\tau_{k-1},\tau_k],\; x\in[-R,R]^r,\\
    \displaystyle \rho(x,t)=0,&  x\in\partial([-R,R]^r).
  \end{cases}
  \]
    with initial condition \(\rho(x,\tau_{k-1})\) and zero Dirichlet boundary conditions on \(\partial([-R,R]^r)\).  This produces the propagated density \(\rho(x,\tau_k)\).  Here, we will design a new PDE solver based on low‐discrepancy sampling for this high‐dimensional Kolmogorov equations; see Section \ref{sec:Kolmogorov-solver} for details.

    Unlike the classical \yauyau algorithm, the prediction step on each interval \([\tau_{k-1},\tau_k]\) is performed by solving the (modified) Kolmogorov forward equation  
\[
  \frac{\partial \rho(x,t)}{\partial t}
  = \mathcal{L}^*\,\rho(x,t)
    - \tfrac12\,\|h(x)\|^2\,\rho(x,t),
\]
which embeds the observation model \(h(x)\) directly into the dynamics.  By contrast, our improved \yauyau scheme solves the standard Kolmogorov equation  
\[
  \frac{\partial \rho(x,t)}{\partial t}
  = \mathcal{L}^*\,\rho(x,t)
\]
omitting the quadratic term \(-\tfrac12\|h(x)\|^2\,\rho\).  This decouples the offline propagation from any observation‐dependent contributions (neither \(h\) nor \(Y\) appear), and confines all observation updates to a simple, multiplicative online correction.  

\paragraph{Online stage (Correction)}
    Upon receiving the observation increment \(\Delta y_k = y_{\tau_k}-y_{\tau_{k-1}}\), update the density by
    \[
      \widetilde\rho_k(x)
      = \exp\!\Bigl\{h(x)^\top \Delta y_k
        - \tfrac12\,\Delta t\,\|h(x)\|^2\Bigr\}\,\rho(x,\tau_k).
    \]
    This update comprises:
    \begin{itemize}
      \item \textbf{Linear term:} \(h(x)^\top\Delta y_k\) adjusts for the observation increment.
      \item \textbf{Quadratic term (Itô correction):} \(-\tfrac12\Delta t\,\|h(x)\|^2\) compensates log‐domain noise bias, which is omitted in the classical \yauyau scheme. Its necessity in the improved \yauyau\ algorithm is explained in Remark~\ref{rem:itocorrection}.
    \end{itemize}

The new \yauyau\ framework introduces three key innovations:
\begin{enumerate}
  \item \textbf{Novel offline-online update:} separates pure‐dynamics propagation from the observation update, simplifying the PDE solve.
  \item \textbf{Low‐discrepancy sampling:} employs Halton or Sobol points for quasi‐uniform state‐space coverage, reducing sampling variance. 
  \item \textbf{Log‐domain correction and acceleration:} includes both linear and quadratic terms for bias-free updates, and precomputes the transition operator offline to enable GPU-accelerated matrix-vector multiplies.
\end{enumerate}
These enhancements preserve the classical prediction-correction structure of the \yauyau\ filter, while supporting scalable high-dimensional implementation through offline low-discrepancy sampling and lightweight online exponential updates. 

\subsubsection{Procedure for the Improved \yauyau\ Algorithm} Now, we summarize the detailed procedure for the \textbf{Improved \yauyau\ Algorithm} as follows:
\begin{enumerate}
  \item \textbf{Initialization:}
    \begin{itemize}
      \item Choose domain $[-R,R]^r$ including the initial value $x_0$ and generate samples $\{x_i\}_{i=1}^n \subset [-R,R]^r$ via a low-discrepancy sequence (Halton/Sobol).
      \item Set $y_{0}=0$, $\rho_0(x)=\sigma_0(x)$, and partition $[0,T]$ into $K$ equal steps $0=\tau_0<\tau_1<\cdots<\tau_K=T$, with $\Delta t=T/K$.
    \end{itemize}

  \item \textbf{Offline stage (Kolmogorov solve):}\\
  For each $k=1,\dots,K$, solve
  \begin{equation}
  \label{eq:Kolmogorov_new}
  \begin{cases}
    \displaystyle \frac{\partial \rho_k(x,t)}{\partial t}
    = \mathcal{L}^*\,\rho_k(x,t), & t\in(\tau_{k-1},\tau_k],\ x\in[-R,R]^r,\\[1mm]
    \rho_k(x,t)=0, & x\in\partial([-R,R]^r),
  \end{cases}
  \end{equation}
  where initial condition is \(\rho(x,\tau_{k-1})\) with $\rho(x,\tau_0) = \sigma_0(x)$ and zero Dirichlet boundary conditions are imposed on $\partial([-R,R]^r)$, to obtain $\rho_k(x,\tau_k)$.

  \item \textbf{Online update:}\\
  Upon receiving observation $y_{\tau_k}$, set
  \[
    \Delta y_k = y_{\tau_k} - y_{\tau_{k-1}},
  \]
  and update
  \[
    \widetilde\rho_k(x)
    = \exp\!\Bigl\{h(x)^\top\Delta y_k
      - \tfrac12\,\Delta t\,\|h(x)\|^2\Bigr\}
      \,\rho_k(x,\tau_k).
  \]
    \item \textbf{Weighting, normalization, and estimation:}\\
After the online update, the unnormalized density \(\widetilde\rho_{k}(x)\) at time \(\tau_k\) represents the posterior \(\rho(x,\tau_k)\). To compute the conditional expectation
\[
  \mathbb{E}[X_{\tau_k}\mid\mathcal{Y}_{\tau_k}]
  = \frac{\displaystyle\int_{[-R,R]^r}x\,\rho(x,\tau_k)\,\mathrm{d}x}
         {\displaystyle\int_{[-R,R]^r}\rho(x,\tau_k)\,\mathrm{d}x},
\]
we use Monte Carlo integration on the sample points \(\{x_i\}_{i=1}^n\):
\[
  \mathbb{E}[X_{\tau_k}\mid\mathcal{Y}_{\tau_k}]
  \approx \frac{\sum_{i=1}^n x_i\,\widetilde\rho_{k}(x_i)}
             {\sum_{i=1}^n \widetilde\rho_{k}(x_i)}.
\]
Equivalently, define the unnormalized and normalized weights
\[
  w_i^{(k)} = \widetilde\rho_{k}(x_i),
  \quad
  \widehat w_i^{(k)} = \frac{w_i^{(k)}}{\sum_{j=1}^n w_j^{(k)}},
\]
and estimate the state by
\[
  \widehat X_{\tau_k}
  = \sum_{i=1}^n \widehat w_i^{(k)}\,x_i.
\]
Here \(w_i^{(k)}\) are the unnormalized weights, \(\widehat w_i^{(k)}\) their normalized counterparts, and the weighted sum \(\widehat X_{\tau_k}\) provides a Monte Carlo approximation to the true conditional expectation \(\mathbb{E}[X_{\tau_k}\mid\mathcal{Y}_{\tau_k}]\).
\end{enumerate}

\section{Convergence Analysis}
\label{sec:convergence-analysis}
Consider two operators
\[
A\sigma = \mathcal{L}^*\,\sigma,
\qquad
B\sigma = \sigma\,h(x)^{\top}\,\dot Y(t).
\]
Denote by \(T_{\Delta t}\) the exact evolution operator over \(\Delta t\):
\[
T_{\Delta t} := \exp\{(A+B)\,\Delta t\},
\]
and our proposed offline-online update reads:
\[
S_{\Delta t}
:= \exp\!\Bigl\{h(x)^{\top}\Delta y
  - \tfrac12\,\Delta t\,\|h(x)\|^2\Bigr\}
  \,\exp\{A\,\Delta t\}.
\]

\begin{theorem}[Local truncation error]\label{thm:local_error}
Under suitable regularity conditions (see Assumption~\ref{assump:regularization}), there exists a constant \(C>0\), independent of \(\Delta t\), such that
\[
\bigl\|T_{\Delta t}\,\sigma(t)\;-\;S_{\Delta t}\,\sigma(t)\bigr\|
\;\le\; C\,\Delta t^2.
\]
That is, the scheme has local truncation error \(O(\Delta t^2)\).
\end{theorem}
The proof of Theorem \ref{thm:local_error} is provided in Appendix~\ref{Appsec:proofthm:local_error}, where we assume the following regularity conditions:
\begin{assumption}[Regularity Conditions]\label{assump:regularization}\leavevmode
\begin{itemize}
  \item The operators \(A\), \(B\), and \(A+B\) generate strongly continuous semigroups on the chosen functional space (e.g.\ \(L^{\infty}\) or a suitable Sobolev space), and the solution \(\sigma(t)\) remains in their domains throughout \([0,T]\).
  \item The initial density \(\sigma_0\) is sufficiently smooth so that the expansion
    \[
      \exp\{(A+B)\Delta t\}
      = I + (A+B)\,\Delta t
      + \tfrac12\,(A+B)^2\,\Delta t^2
      + O(\Delta t^3)
    \]
    holds in the chosen operator norm.
  \item The observation function \(h(x)\) is continuously differentiable and bounded, ensuring that the exponential update factor
    \(\exp\{h(x)^\top\Delta y - \tfrac12\Delta t\,\|h(x)\|^2\}\)
    is well defined and smooth in \(x\).
  \item The commutator \([B,A] := BA - AB\) and higher‐order nested commutators are bounded in the same norm, so that the \(O(\Delta t^3)\) remainder can be uniformly controlled.
\end{itemize}
\end{assumption}

\begin{theorem}[Global convergence]\label{thm:convergence}
Let \(\sigma(T)\) be the exact solution of the DMZ \eqref{eq:DMZ} at time \(T\) and \(\sigma^\Delta(T)\) the result of applying \(K=T/\Delta t\) steps of the operator \(S_{\Delta t}\).  Under the regularity conditions (Assumption~\ref{assump:regularization}) and the stability conditions (Assumption~\ref{assump:stability}), there exists \(C_1>0\) such that
\[
\bigl\|\sigma(T) - \sigma^\Delta(T)\bigr\|
\;\le\; C_1\,\Delta t.
\]
Hence the method is first‐order accurate globally.
\end{theorem}

The proof of \Cref{thm:convergence} is deferred to Appendix~\ref{Appsec:proofthm:convergence}, where we require that:
\begin{assumption}[Stability Conditions]\label{assump:stability}
The operator $S_{\Delta t}$ satisfies
\[
\|S_{\Delta t}\|\;\le\;1 + L\,\Delta t
\]
for some constant \(L>0\). This ensures errors do not grow uncontrollably.
\end{assumption}

\begin{remark}[Necessity of the Itô correction]\label{rem:itocorrection}
Consider the naive update factor \(\exp\{h(x)^{\top}\Delta y\}\) with
\(\Delta y\sim\mathcal{N}(0,\Delta t\,I)\).  Let \(Z = h(x)^{\top}\Delta y\), then $Z$ follows the Gaussian distribution with zero mean and \[
\operatorname{Var}(Z)=\Delta t\,\|h(x)\|^2.
\]
Hence, \[\mathbb{E}[e^Z]=\exp\{\tfrac12\Delta t\,\|h(x)\|^2\}\neq1,\]
so using \(\exp\{h(x)^{\top}\Delta y\}\) alone introduces a systematic scaling bias.  
By including the correction term \(-\tfrac12\Delta t\,\|h(x)\|^2\), the updated factor
\(\exp\{h(x)^{\top}\Delta y - \tfrac12\Delta t\,\|h(x)\|^2\}\)
has expectation one, matching the exact Itô expansion and ensuring an unbiased update.
\end{remark}

We conclude from Theorems~\ref{thm:local_error} and \ref{thm:convergence} that the proposed new \yauyau scheme achieves a local truncation error of order \(O(\Delta t^2)\) and a global error of order \(O(\Delta t)\), i.e.\ first‐order accuracy in the time step. Under the regularity and stability assumptions, the operators \(A\) and \(B\) generate well‐posed semigroups, and the inclusion of the Itô correction term guarantees an unbiased, statistically consistent update in accordance with the DMZ equation.

\section{Solving Kolmogorov Equations via Low-Discrepancy Sampling and Kernel Approximation}\label{sec:Kolmogorov-solver}

Consider the Kolmogorov forward equation
\begin{equation}\label{eq:Kolmogorov}
  \frac{\partial \rho(x,t)}{\partial t}
  = \mathcal{L}^*\,\rho(x,t),
\end{equation}
where the adjoint operator is defined by
\[
  \mathcal{L}^*\,\rho
  = \frac{1}{2}\,\Delta\rho
    - \nabla\cdot \bigl(f(x)\,\rho\bigr).
\]
Classical grid‐based solvers suffer from the \emph{curse of dimensionality}, with complexity growing exponentially in the state dimension $r$ (see e.g., \cite{YuehLinYau2014}). To overcome this, we employ Halton- or Sobol-sequence low-discrepancy sampling and propose a novel \textbf{kernel-based operator approximation}, leveraging the uniform coverage of the sample points to  solve \eqref{eq:Kolmogorov} via efficient matrix operations.

\subsection{Discrete Solution-Mapping Operator Construction}
We seek a solution‐mapping operator \(L_{\Delta t}\) (cf. time‐evolution operator, Fokker–Planck propagator, or transition operator in the literature) such that, for any sufficiently smooth test function $u$,
\[
  \frac{L_{\Delta t}u - u}{\Delta t}
  \;\longrightarrow\;
  \mathcal{L}^*\,u
  \quad\text{as }\Delta t\to0.
\]
Equivalently, $L_{\Delta t}$ propagates the initial condition $u$ forward by one time step $\Delta t$.  

Using a kernel‐based approach, we discretize $L_{\Delta t}$ on $n$ sample points $\{x_i\}_{i=1}^n$.  First, we approximate $u$ in the basis of Dirac masses $\{\delta_{x_i}(x)\}_{i=1}^n$:
\[
  u(x) \;\approx\; \sum_{i=1}^n a_i\,\delta_{x_i}(x),
\]
where \[
  \delta_{x_i}(x) \;=\; \delta\bigl(x - x_i\bigr),\quad a_i = u(x_i),
  \quad i=1,\dots,n.
\] Then defining $$f_{\Delta t,i} := L_{\Delta t}\,\delta_{x_i},$$ whose approximation in the Dirac basis is 
\[
  f_{\Delta t,i}
  \;\approx\;
  \sum_{j=1}^n f_{\Delta t,i}(x_j)\,\delta_{x_j},
\]
and hence the linearity of the Kolmogorov equation gives
\begin{align*}
L_{\Delta t} u &\approx L_{\Delta t} \sum_{i=1}^n a_i \delta_{x_i}  = \sum_{i=1}^n a_i L_{\Delta t} \delta_{x_i}
= \sum_{i=1}^n a_i f_{{\Delta t}, i} \approx \sum_{i=1}^n a_i \sum_{j=1}^n f_{{\Delta t}, i}(x_j) \delta_{x_j}.    
\end{align*}
Therefore,
\begin{equation}
    \label{eq:L_Delta_t}
  \boxed{L_{\Delta t}u
  \;\approx\;
  \sum_{j=1}^n \Bigl(\sum_{i=1}^n a_i\,f_{\Delta t,i}(x_j)\Bigr)\,\delta_{x_j}.}
\end{equation}
Defining the coefficient vector $a=[a_1,\ldots,a_n]^{\top}$, \Cref{eq:L_Delta_t} corresponds to a linear mapping
\[
  a\;\mapsto\;F_{\Delta t}^{\!\top}\,a,
\]
where the $(i,j)$ element of the matrix $F_{\Delta t}$ is defined by:
\[\boxed{F_{\Delta t}(i,j) = f_{\Delta t,i}(x_j).}\]
We can approximate $f_{\Delta t,i}(x_j)$ via some continuous kernel $K_{\Delta t}(x,y)$. Then the continuous operator is
\[
  (L_{\Delta t}u)(x)
  = \int_{\mathbb{R}^r}K_{\Delta t}(x,y)\,u(y)\,dy,
\]
and one can show that
\[
  \frac{L_{\Delta t}u(x)-u(x)}{\Delta t}
  = \mathcal{L}^*u(x)
  + O(\sqrt{\Delta t}),
\]
as $\Delta t\to0$ (see \Cref{thm:1storderkernel}), where 
\[K_{\Delta t}(x,y)
  = c_{\Delta t}
    \exp\!\Bigl(
      -\tfrac{\|x-y\|^2}{2\Delta t}
      - (y-x)\cdot f(x) - \Delta t\bigl(\nabla\cdot f(x) + \tfrac12\|f(x)\|^2\bigr)
    \Bigr),\]
with normalization coefficient
\[
  c_{\Delta t}
  = \frac{1}{(2\pi\,\Delta t)^{r/2}}.
\]
A higher-order kernel function can be also established which will be discussed later in Section \ref{sec:high-order-kernel}.

\begin{remark}
In $K_{\Delta t}(x,y)$, the first argument $x$ is the evaluation point, and $y$ is the integration variable (a sample point).  On the discrete sample set $\{x_i\}$, the integral reduces to
\[
  (L_{\Delta t}u)(x_i)
  \approx
  \sum_{j=1}^n K_{\Delta t}(x_i,x_j)\,u(x_j)\,\Delta y,
\]
so that
\[
  \boxed{f_{\Delta t,j}(x_i)
  = K_{\Delta t}(x_i,x_j),\quad \forall i,j\in \{1,\ldots,n\}.}
\]
\end{remark}

\subsection{Algorithm for Kolmogorov Equation}

The main steps are:
\begin{enumerate}
  \item \textbf{Sampling:} Generate $n$ quasi‐uniform points $\{x_i\}_{i=1}^n$ in $[-R,R]^r$ using a Sobol or Halton low-discrepancy sequence.
  \item \textbf{Operator matrix assembly:} For each pair $(i,j)\in \{1,\ldots,n\}^2$, compute
    \[
      F_{\Delta t}(i,j)
      = K_{\Delta t}(x_j,x_i),
    \]
    forming the discrete solution-mapping operator $F_{\Delta t}$.
  \item \textbf{Boundary conditions:} Impose zero Dirichlet conditions on $\partial([-R,R]^r)$ by identifying any sample point within tolerance $\omega$ (e.g., $\omega = 10^{-12}$) of the boundary and zeroing out the corresponding row of $F_{\Delta t}$.
  \item \textbf{Initial condition discretization:} Evaluate the initial density $\sigma_0(x)$ at the sample points to form
    \[
      \sigma^{(0)}
      = \bigl[\sigma_0(x_1),\dots,\sigma_0(x_n)\bigr]^\top.
    \]
  \item \textbf{Time stepping via matrix multiplication:}
    \[
      \sigma^{(k+1)} = F_{\Delta t}^{\!\top}\,\sigma^{(k)},
      \quad k = 0,1,\dots.
    \]
\end{enumerate}

This approach yields a scalable, kernel‐based solver for high‐dimensional Kolmogorov equations, exploiting the quasi‐uniform, low‐variance coverage of low‐discrepancy samples alongside efficient linear‐algebra operations.

\subsection{Kernel Approximation}\label{sec:high-order-kernel}

In this section, we construct the kernel function to approximate the solution-mapping operator of the Kolmogorov forward equation. 

\subsubsection{First-Order Kernel Construction}
Let us start with a first-order approximation of the adjoint Kolmogorov operator $\mathcal{L}^*$ over a time step size $\Delta t$, resulting in a local truncation error of \(O(\sqrt{\Delta t})\).
\begin{theorem}[First-Order Approximation of the Adjoint Kolmogorov Operator]\label{thm:1storderkernel}
Let the continuous integral operator be defined by
\[
  (L_{\Delta t}u)(x)
  = \int_{\mathbb{R}^r}
      K_{\Delta t}(x,y)\,u(y)\,dy,
\]
where the kernel function is given by
\begin{equation}
  K_{\Delta t}(x,y)
  = c_{\Delta t}
    \exp\!\Bigl(
      -\tfrac{\|x-y\|^2}{2\Delta t}
      - (y-x)\cdot f(x)
      - \tfrac{\Delta t}{2}\bigl(2\,\nabla\cdot f(x)+\|f(x)\|^2\bigr)
    \Bigr), \label{eq:kernel}
\end{equation}
with normalization constant
\[
  c_{\Delta t}
  = \frac{1}{(2\pi\,\Delta t)^{r/2}},
  \quad
  \int_{\mathbb{R}^r}K_{\Delta t}(x,y)\,\mathrm{d}y= 1.
\]
Let the adjoint Kolmogorov operator be
\[
  \mathcal{L}^*u(x)
  = \tfrac12\,\Delta u(x)
    - \nabla\cdot \bigl(f(x)\,u(x)\bigr).
\]
Then, for any sufficiently smooth \(u\), as \(\Delta t\to0\),
\[
  \frac{(L_{\Delta t}u)(x) - u(x)}{\Delta t}
  = \mathcal{L}^*u(x)
    + O\bigl(\sqrt{\Delta t}\,\bigr).
\]
\end{theorem}
The proof of Theorem \ref{thm:1storderkernel} is provided in Appendix~\ref{Appsec:proofthm:1storderkernel}.

\subsubsection{Second-Order Kernel Construction}
We show that by taking suitable linear combinations of the first‐order kernel evaluated at different time steps together with the Dirac operator, one can construct a second order kernel with an improved local truncation error of \(O(\Delta t)\).

\begin{theorem}[Multi‐Time‐Scale Kernel: Second‐Order Consistency]\label{thm:2ndorderkernel}
Let the kernel \(K_{\Delta t}(x,y)\) and its integral operator
\[
  (L_{\Delta t}u)(x)
  = \int_{\mathbb{R}^r}K_{\Delta t}(x,y)\,u(y)\,\mathrm{d}y
\]
admit the first‐order expansion as described in \Cref{thm:1storderkernel}: 
\[
  (L_{\Delta t}u)(x)
  = u(x) + \Delta t\,\mathcal{L}^*u(x)
    + C\,(\Delta t)^{3/2} + O(\Delta t^2),
\]
with \(C\neq0\).  Define the multi‐time‐scale operator
\[
  (L^{(2)}_{\Delta t}u)(x)
  = \alpha\,L_{\Delta t}u(x)
    + \beta\,L_{\Delta t/2}u(x)
    + \gamma\,u(x),
\]
with coefficients
\[
  \alpha = -\sqrt{2}-1,\quad
  \beta  = 4+2\sqrt{2},\quad
  \gamma = -2-\sqrt{2}.
\]
Then, as \(\Delta t\to0\),
\[
  \frac{(L^{(2)}_{\Delta t}u)(x)-u(x)}{\Delta t}
  = \mathcal{L}^*u(x) + O(\Delta t),
\]
Consequently, the corresponding kernel
\[
  K^{(2)}_{\Delta t}(x,y)
  = \alpha\,K_{\Delta t}(x,y)
    + \beta\,K_{\Delta t/2}(x,y)
    + \gamma\,\delta(x-y)
\]
approximates the adjoint Kolmogorov operator to second order.
\end{theorem}
The proof is presented in Appendix~\ref{Appsec:proofthm:2ndorderkernel}.

\subsubsection{High-Order Kernel Construction}
The idea developed in Theorem \ref{thm:2ndorderkernel} can be further extended to arbitrary high-order. Let \(s_1,\dots,s_N\) be distinct positive scale factors, e.g., \[s_1=1,\;s_2=2,\;s_3=4,\;\dots,\;s_N=2^{N-1}.\]  By \Cref{thm:1storderkernel}, for each \(j=1,\dots,N\), let
\[
  (L_{\Delta t/s_j}u)(x)
  = u(x)
    + \frac{\Delta t}{s_j}\,\mathcal{L}^*u(x)
    + C_1\Bigl(\tfrac{\Delta t}{s_j}\Bigr)^{\mu_1}
    + \cdots
    + C_{N-1}\Bigl(\tfrac{\Delta t}{s_j}\Bigr)^{\mu_{N-1}}
    + O\bigl((\Delta t)^{\mu_N}\bigr),  
\]
where each \(C_k\neq0\), and exponents \(\mu_1<\mu_2<\cdots<\mu_{N-1}\) denote the lower‐order error terms to be eliminated (for instance, \(\mu_1=\tfrac32,\;\mu_2=2,\dots\)). To match both the constant term and the \(\Delta t\) coefficient while eliminating the next \(N-1\) error orders, define the operator
\[
  (L^{(N)}_{\Delta t}u)(x)
  = \sum_{j=1}^N c_j\,\bigl(L_{\Delta t/s_j}u\bigr)(x)
    \;+\; c_{N+1}\,u(x),
\]
where the first \(N\) terms are the scaled integral operators and the last term is the identity operator (with kernel Dirac \(\delta(x-y)\) satisfying $u(x) = \int \delta(x-y)u(y)\mathrm{d}y$).  

Requiring the expansion
\[
  (L^{(N)}_{\Delta t}u)(x)
  = u(x)
    + \Delta t\,\mathcal{L}^*u(x)
    + O\bigl((\Delta t)^{\mu_N}\bigr)
\]
imposes the linear constraints on the coefficient vector \(c=[c_1,\dots,c_{N+1}]^\top\):
\begin{enumerate}
  \item \emph{Constant‐term matching:}\quad \(\displaystyle \sum_{j=1}^{N+1}c_j = 1.\)
  \item \emph{\(\Delta t\)‐term matching:}\quad \(\displaystyle \sum_{j=1}^N c_j\,s_j^{-1} = 1.\)
  \item \emph{Error‐term elimination:}\quad for each \(k=1,\dots,N-1\),
    \[\displaystyle \sum_{j=1}^N c_j\,s_j^{-\mu_k} = 0.\]
\end{enumerate}
These \(N+1\) linear equations can be written in matrix form
\[
  V\,c = d,
\]
with
\[
  V = 
  \begin{bmatrix}
    1              & 1              & \cdots & 1              & 1 \\[4pt]
    s_1^{-1}       & s_2^{-1}       & \cdots & s_N^{-1}       & 0 \\[4pt]
    s_1^{-\mu_1}   & s_2^{-\mu_1}   & \cdots & s_N^{-\mu_1}   & 0 \\[2pt]
    \vdots         & \vdots         & \ddots & \vdots         & \vdots \\
    s_1^{-\mu_{N-1}} & s_2^{-\mu_{N-1}} & \cdots & s_N^{-\mu_{N-1}} & 0
  \end{bmatrix},
  \quad
  d =
  \begin{bmatrix}
    1 \\[3pt]
    1 \\[3pt]
    0 \\[2pt]
    \vdots \\[2pt]
    0
  \end{bmatrix}.
\]
Since all \(s_j\) are distinct and nonzero, \(V\) (a generalized Vandermonde form) is nonsingular, and the unique solution is given by Cramer's rule:
\[
  c_j = \frac{\det V_j}{\det V},
  \quad j=1,2,\dots,N+1,
\]
where \(V_j\) is formed by replacing the \(j\)th column of \(V\) with \(d\).

Thus one obtains the high-order operator
\[
  \boxed{(L^{(N)}_{\Delta t}u)(x)
  = \sum_{j=1}^N c_j\,(L_{\Delta t/s_j}u)(x)
    \;+\; c_{N+1}\,u(x),}
\]
and the corresponding high‐order kernel
\begin{equation}
    \label{eq:general_kernel}
    \boxed{
    K^{(N)}_{\Delta t}(x,y)
    = \sum_{j=1}^N c_j\,K_{\Delta t/s_j}(x,y)
      \;+\; c_{N+1}\,\delta(x-y).
  }
\end{equation}
Obviously, this kernel is general for any high-order $N$ and recovers the first‐order kernel of \Cref{thm:1storderkernel} (with \(N=1\)) and the second‐order kernel of \Cref{thm:2ndorderkernel} (with \(N=2\)).

\section{Error Analysis under Low‐Discrepancy Sampling}\label{sec:error_analysis_qmc}

In Section~\ref{sec:convergence-analysis}, we established the error bounds and convergence of new \yauyau scheme under the assumption that the evolution operators \(\exp\{A\Delta t\}\) and \(\exp\{B\Delta t\}\) are available exactly.  In the improved \yauyau algorithm, however, these operators are replaced by discrete approximations evaluated on low‐discrepancy samples \(\{x_i\}_{i=1}^n\subset[-R,R]^r\).  It is therefore essential to quantify the additional error introduced by this QMC discretization and to show that the overall method preserves similar convergence properties.  In this section, we derive the corresponding QMC‐based error bounds and convergence estimates.

\begin{definition}[Hardy–Krause variation]
Let \(\Omega=[-R,R]^r\). For any nonempty index set \(\mathcal{U}=\{i_1,\dots,i_k\}\subseteq\{1,\dots,r\}\),
define the restriction
\[
g_{_{\mathcal{U}}}(x_{i_1},\dots,x_{i_k})
= g(x_1,\dots,x_r)\bigl|_{x_j=R\,(j\notin \mathcal{U})},
\]
and its Vitali variation
\[
V(g_{_{\mathcal{U}}})
= \int_{[-R,R]^k}
\Bigl|\frac{\partial^k g_{_{\mathcal{U}}}}{\partial x_{i_1}\cdots\partial x_{i_k}}\Bigr|
\,\mathrm{d}x_{i_1}\cdots\mathrm{d}x_{i_k}.
\]
The Hardy–Krause variation is then
\[
V_{\mathrm{HK}}(g)
= \sum_{\emptyset\neq \mathcal{U}\subseteq\{1,\dots,r\}}V(g_{_{\mathcal{U}}}).
\]
\end{definition}

\begin{lemma}[Koksma–Hlawka inequality {\cite{Niederreiter1992}}]
Let \(g:\Omega\to\mathbb{R}\) have Hardy–Krause variation \(V_{\mathrm{HK}}(g)\).  Sample
\(\{x_i\}_{i=1}^n\subset\Omega\) by a low‐discrepancy sequence with
\[
D^*(n)=O\!\Bigl(\frac{(\log n)^r}{n}\Bigr).
\]
Then
\[
\Bigl|\frac1n\sum_{i=1}^n g(x_i)
-\frac1{(2R)^r}\int_\Omega g(x)\,\mathrm{d}x\Bigr|
\;\le\; V_{\mathrm{HK}}(g)\,D^*(n).
\]
\end{lemma}

\subsection{Offline QMC Error (Approximation of \(\exp\{A\Delta t\}\))}

\begin{definition}[Operator norm \(\|\cdot\|_{\infty\to\infty}\)]
For \(T:L^\infty(\Omega)\to L^\infty(\Omega)\), define
\[
\|T\|_{\infty\to\infty}
=\sup_{\|\phi\|_\infty\le1}\|T\phi\|_\infty,
\quad
\|\phi\|_\infty=\sup_{x\in\Omega}|\phi(x)|.
\]
\end{definition}
\begin{remark}
We choose the \(L^\infty\to L^\infty\) operator norm because the Koksma–Hlawka and related QMC error bounds are naturally formulated in terms of supremum (worst‐case) norms, allowing direct control of the maximum approximation error over the entire domain.
\end{remark}
\begin{definition}[Exact offline operator]
For \(\phi\in L^\infty(\Omega)\), set
\[
(\exp\{A \Delta t\}\phi)(y)
=\frac1{(2R)^r}\int_\Omega K_{\Delta t}(x,y)\,\phi(x)\,\mathrm{d}x.
\]
\end{definition}

\begin{definition}[QMC‐approximate offline operator]
With the same low‐discrepancy samples \(\{x_i\}_{i=1}^n\subset\Omega\),
define
\[
(\widetilde T_n\phi)(y)
=\frac1n\sum_{i=1}^n K_{\Delta t}(x_i,y)\,\phi(x_i).
\]
\end{definition}

\begin{theorem}[Offline QMC error]\label{thm:offline_qmc}
Assume for each fixed \(y\), the map \(x\mapsto K_{\Delta t}(x,y)\) has bounded
Hardy–Krause variation and set
\[
C_A=\sup_{y}V_{\mathrm{HK}}(K_{\Delta t}(\cdot,y))<\infty.
\]
Then
\[
\|\exp\{A \Delta t\}-\widetilde T_n\|_{\infty\to\infty}
\;\le\; C_A\,D^*(n).
\]
\end{theorem}
The proof is deferred to Appendix~\ref{Appsec:proofthm:offline_qmc}.

\subsection{Online QMC Error (Approximation of \(\exp\{B\Delta t\}\))}

\begin{definition}[Exact observation update operator]
Let \(\Omega=[-R,R]^r\). For any test function \(\phi\in L^\infty(\Omega)\), define
\[
g_y(x)
=\exp\!\{h(x)^\top\Delta Y-\tfrac12\Delta t\,\|h(x)\|^2\},
\] and the exact observation update operator is 
\[
\bigl(\exp\{B\Delta t\}\,\phi\bigr)(y)
=\frac{1}{(2R)^r}
\int_{\Omega}
g_y(x)
\,\phi(x)\,\mathrm{d}x.
\]
\end{definition}

\begin{definition}[QMC approximate update operator]
Using the same low‐discrepancy samples \(\{x_i\}_{i=1}^n\subset\Omega\) with star-discrepancy
\[
D^*(n)=O\!\bigl((\log n)^r/n\bigr),
\]
define the QMC approximate update operator
\[
\bigl(\widetilde S_n\phi\bigr)(y)
=\frac{1}{n}
\sum_{i=1}^n
\exp\!\bigl(h(x_i)^{\top}\,\Delta Y - \tfrac12\Delta t\,\|h(x_i)\|^2\bigr)
\,\phi(x_i).
\]
\end{definition}

\begin{theorem}[Online QMC error]\label{thm:online_qmc}
Assume for each fixed \(y\), the map \(x\mapsto g_y(x)\) has bounded
Hardy–Krause variation \(V_{\mathrm{HK}}(g_y)\le C_B<\infty\).  Then the exact
measurement update and its QMC approximation satisfy
\[
\bigl\|\exp\{B\Delta t\}-\widetilde S_n\bigr\|_{\infty\to\infty}
\;\le\;
C_B\,D^*(n).
\]
\end{theorem}
The proof is deferred to Appendix~\ref{Appsec:proofthm:online_qmc}.
\subsection{Local and Global QMC Error Estimates}
\begin{theorem}[Local error estimate]\label{thm:qmc_local}
Under the hypotheses of
Theorems~\ref{thm:local_error}, \ref{thm:offline_qmc} and \ref{thm:online_qmc},
there exist constants \(C_1,C_2>0\) such that
\[
\|T_{\Delta t}\sigma - S_{\Delta t}^n\sigma\|
\;\le\;
C_1\,\Delta t^2 + C_2\,D^*(n).
\]
\end{theorem}
The proof is deferred to Appendix~\ref{Appsec:proofthm:qmc_local}.

\begin{theorem}[Global error estimate]\label{thm:qmc_global}
Under the local truncation error bound of \Cref{thm:qmc_local} and the stability
assumption \(\|S_{\Delta t}^n\|\le1+L\Delta t\), there exists \(C_3>0\) such that
for \(T=K\Delta t\),
\[
\|\sigma(T)-(S_{\Delta t}^n)^K\sigma(0)\|
\;\le\;
C_3\Bigl(\Delta t + \frac{D^*(n)}{\Delta t}\Bigr).
\]
\end{theorem}
The proof is deferred to Appendix~\ref{Appsec:proofthm:qmc_global}.
    
\begin{figure}[ht]
  \centering
  \includegraphics[width=0.48\textwidth]{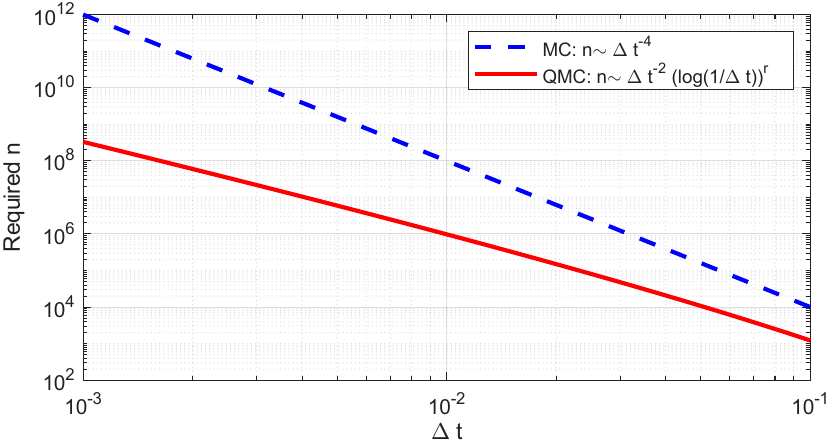}
  \caption{Required sample size \(n\) vs.\ \(\Delta t\) on a log-log scale,
    comparing MC (\(n\sim\Delta t^{-4}\), dashed blue) and 
    QMC (\(n\sim\Delta t^{-2}(\log(1/\Delta t))^3\), solid red) for \(r=3\).}
  \label{fig:n_vs_dt_mc_qmc}
\end{figure}

In the worst‐case scenario, one balances the sampling discrepancy \(D^*(n)\) with the error \(\Delta t^2\) to preserve a local truncation error of order \(O(\Delta t^2)\) and a global error of order \(O(\Delta t)\).  For standard Monte Carlo (MC) with \(D^*(n)=O(n^{-1/2})\), this leads to
\[
n_{\rm MC}\sim \Delta t^{-4},
\]
whereas for quasi‐Monte Carlo (QMC) with \(D^*(n)=O\!\bigl((\log n)^r/n\bigr)\) one obtains
\[
n_{\rm QMC}\sim \Delta t^{-2}\,\bigl(\log(1/\Delta t)\bigr)^r.
\]
Thus, for a given dimension \(r\), QMC requires asymptotically far fewer samples than MC to achieve the same accuracy. Fig.~\ref{fig:n_vs_dt_mc_qmc} illustrates this comparison for \(r=3\). Although the asymptotic estimates predict extremely large sample sizes in the worst case, in practice the number of samples needed to achieve a desired error is often orders of magnitude smaller.

\section{Local Resampling-Restart Mechanism}
\label{sec:local-resampling-restart}

\begin{figure*}[ht]
  \centering
  \includegraphics[width=\textwidth]{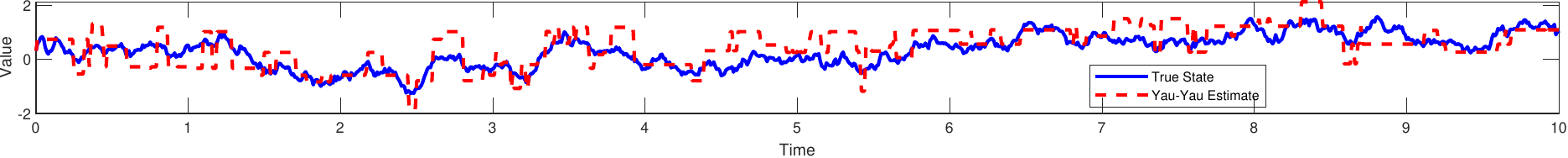}
  \caption{The “great‐wall” phenomenon: due to excessively sparse sampling, the filter correction stagnates along a coordinate direction.}
  \label{fig:greatwall}
\end{figure*}

In high‐dimensional nonlinear filtering, globally distributed samples can leave large regions of the state space severely under‐covered, giving rise to the so‐called “great‐wall” phenomenon: within a single time step \(\Delta t\), the lack of local sample support causes the filter update to stagnate along one or more coordinate directions (see Fig.~\ref{fig:greatwall}).

To overcome this, we propose a \emph{local resampling-restart} strategy into our improved \yauyau algorithm:

\paragraph{Local resampling} After running the improved \yauyau filter for one or more time steps to obtain the current state estimate \(\widehat x_k\), we recenter a hypercube
\(\bigl[\widehat x_k - R,\;\widehat x_k + R\bigr]^r\) around \(\widehat x_k\) and within this local region, we generate \(n_R\) quasi-uniform samples using a low-discrepancy sequence.
The half-width \(R\) and sample count \(n_R\) can be set a priori or adjusted adaptively to ensure sufficient local density.
To optimize computational overhead, one can pre-generate a single reference set
\(\{x_i\}_{i=1}^{n_R}\subset[-R,R]^r\)
and translate it by \(\widehat x_k\) at each restart.
This strategy preserves the low-discrepancy structure while avoiding repeated sampling.

\paragraph{Filter restart}  
After obtaining the new estimate \(\widehat x_k\), we reinitialize the improved \yauyau algorithm within the localized hypercube
\(
\bigl[\widehat x_k - R,\;\widehat x_k + R\bigr]^r
\)
using \(\widehat x_k\) as the new starting point. 

Through this local resampling-restart strategy, the “great-wall” artifact is effectively eliminated, yielding substantial gains in both accuracy and robustness for high-dimensional problems. The rationale is as follows: suppose \(n_G\) samples are drawn uniformly over the global hypercube \(\bigl[-R_G,R_G\bigr]^r\). By uniformity, the expected count \(n_R\) of points inside a smaller local hypercube \(\bigl[-R,R\bigr]^r\) satisfies
\begin{equation}
\label{eq:ruleoflocalvsglobalregion}
\frac{n_R}{n_G}
=
\frac{\operatorname{Vol}\bigl([-R,R]^r\bigr)}{\operatorname{Vol}\bigl([-R_G,R_G]^r\bigr)}
=
\Bigl(\frac{R}{R_G}\Bigr)^{r}.
\end{equation}
Maintaining a fixed local sample size \(n_R\) via purely global sampling therefore requires
\[
n_G
=
n_R \,\Bigl(\frac{R_G}{R}\Bigr)^{r},
\]
which grows exponentially in \(r\). By contrast, directly resampling \(n_R\) points within \(\bigl[-R,R\bigr]^r\) concentrates samples where needed, avoiding this exponential cost. For example, with \(n_R=100\), \(R=0.1\), \(R_G=1\), and \(r=10\),
\[
n_G = 100 \times \Bigl(\tfrac{1}{0.1}\Bigr)^{10} = 10^{12}.
\]
Hence, 100 local samples are equivalent to \(10^{12}\) global samples in terms of local coverage for problems of dimension $10$, demonstrating the effectiveness of the local resampling-restart strategy in overcoming the curse of dimensionality.

\begin{remark}
  The convergence of the improved \yauyau scheme with local resampling-restart can be derived straightforwardly by partitioning the time interval \([0,T]\) into \(K = T/T_0\) segments of equal length \(T_0 = \iota\,\Delta t\) for some restart interval $\iota\geq 1$.  If $\iota$ is not too large, then on each segment, the scheme alone achieves a local truncation error of \(O(\Delta t^2)\) (\Cref{thm:local_error}) and a global error of \(O(\Delta t)\) (\Cref{thm:convergence}). The choice of segment length \(T_0\) impacts practical performance:
  \begin{itemize}
    \item If \(T_0 = \Delta t\) (i.e., $\iota=1$), the algorithm reduces to the original version without restart, yielding global error \(O(\Delta t)\).
    \item For moderate \(\iota\) (e.g.\ \(2\le \iota\le4\)), per‐segment error remains small, preserving the \(O(\Delta t)\) convergence with few restarts.
    \item Very large \(\iota\) allows more error to accumulate within each segment, which can degrade practical accuracy despite fewer restarts.
  \end{itemize}
  Therefore, a practical guideline is: for high-dimensional problems (e.g.\ \(r\ge5\)), balance restart frequency against segment-wise error accumulation by choosing
\[
  T_0 = \iota\,\Delta t,\quad \iota\in[2,5],
\]
which often preserves first-order global convergence while preventing “great-wall” stagnation. Conversely, for low-dimensional problems (e.g.\ \(r<5\)), a single global sampling over a sufficiently large domain \([-R,R]^r\) typically achieves the desired accuracy on \([0,T]\).
\end{remark}

\section{General Nonlinear Filtering with Noise Coefficients}
\label{sec:noise-coeff-filtering}

Consider the more general nonlinear filtering problem, where both the state and observation noises have state‐dependent coefficients:
\begin{equation}
    \label{eq:NonlinearFilter_noise}
\begin{cases}
  \mathrm{d}X_t = f(X_t)\,\mathrm{d}t\;+\; U(X_t)\,\mathrm{d}V_t, 
  & X_t\in\mathbb{R}^r, X_0 = x_0,\\[1ex]
  \mathrm{d}Y_t = h(X_t)\,\mathrm{d}t\;+\; V(X_t)\,\mathrm{d}W_t,
  & Y_t\in\mathbb{R}^q, Y_0 =0
\end{cases}
\end{equation}
with \(U(x)\in \mathbb{R}^{r\times r}\) the state‐noise coefficient and \(V(x)\in \mathbb{R}^{q\times q}\) the observation‐noise coefficient. Let  
\[a(x)=U(x)U(x)^\top,\quad b(x)=V(x)V(x)^\top.\]
We can extend the improved \yauyau algorithm described in Sections~\ref{sec:improved-Yau–Yau-algorithm} (the basic framework) and \ref{sec:local-resampling-restart} (with resampling-restart) to the general nonlinear filtering problem \eqref{eq:NonlinearFilter_noise}, where all modifications (for both invertible and singular cases in $a(x)$ and $b(x)$ respectively) are summarized in Table \ref{tab:generalformulas}. Here, we define 
\[\Delta\bigl(a\sigma\bigr)=\sum_{i,j}\partial^2_{x_i x_j}(a_{ij}\,\sigma).\]
The proofs of these formulas are omitted as they follow the similar arguments as in the case \eqref{eq:NonlinearFilter}.

Note that for singular cases (i.e., either \(a(x)\) or \(b(x)\) is singular), we replace its inverse by the Moore–Penrose pseudoinverse. For example, suppose that $a(x)$ is singular and has an SVD decomposition:
\[
  a(x) = U_a(x)\,\Lambda_a(x)\,U_a(x)^\top,
  \qquad
  \Lambda_a=\mathrm{diag}(\lambda_1,\dots,\lambda_r),
\]
with \(k\) ($< r$) positive eigenvalues \(\{\lambda_i>0\}_{i=1}^k\).  Define the pseudoinverse of \(\Lambda_a\) by
\[
  (\Lambda_a^\dagger)_{ii}
  = \begin{cases}
      1/\lambda_i, & \lambda_i>0,\\
      0,           & \lambda_i=0,
    \end{cases}
\]
then the pseudoinverse of $a(x)$ is given by:
$$
  a(x)^\dagger
  = U_a(x)\,\Lambda_a^\dagger\,U_a(x)^\top.$$
In the kernel, we replace $\|u\|^2_{a(x)^{-1}}$ by \[\|u\|^2_{a(x)^\dagger}=u^\top a(x)^\dagger u\] and
$\det a(x)$ by \[\det\nolimits_{+}a(x)=\prod_{\lambda_i>0}\lambda_i\] (the product of all positive eigenvalues) accordingly. 

\begin{table*}[h]
\centering
\caption{Key differences in the improved Yau–Yau algorithm for the general nonlinear filtering problem \eqref{eq:NonlinearFilter_noise} under invertible versus singular noise‐coefficient matrices}\label{tab:generalformulas}
\resizebox{\textwidth}{!}{%
\renewcommand{\arraystretch}{1.4}
\begin{tabular}{l|p{0.49\textwidth}|p{0.48\textwidth}}
\toprule
 & \textbf{Invertible case} & \textbf{Singular case} \\ \midrule
\textbf{Adjoint operator $\mathcal{L}^*$} &
$\displaystyle
\mathcal{L}^*\sigma
=-\nabla\!\cdot(f\,\sigma)
+\tfrac12\,\Delta(a\,\sigma)
$ &
$\displaystyle
\mathcal{L}^*\sigma
=-\nabla\!\cdot(f\,\sigma)
+\tfrac12\,\Delta(a^\dagger\,\sigma)
$ 
\\[1ex] \hline
\textbf{DMZ} &
$\displaystyle
\mathrm{d}\sigma
=\mathcal{L}^*\sigma\,\mathrm{d}t
+\sigma\,[\,b^{-1}h\,]^\top\,\mathrm{d}Y_t
$ &
$\displaystyle
\mathrm{d}\sigma
=\mathcal{L}^*\sigma\,\mathrm{d}t
+\sigma\,[\,b^\dagger h\,]^\top\,\mathrm{d}Y_t
$ 
\\[1ex] \hline
\textbf{Offline-Online Operators} &
$\displaystyle
A\sigma=\mathcal{L}^*\sigma,\quad
B\sigma=\sigma\,[\,b^{-1}h\,]^\top\,\dot Y(t)
$ &
$\displaystyle
A\sigma=\mathcal{L}^*\sigma,\quad
B\sigma=\sigma\,[\,b^\dagger h\,]^\top\,\dot Y(t)
$ 
\\[2ex] \hline
\textbf{Kernel $K_{\Delta t}(x,y)$} &
\(\displaystyle
\begin{aligned}
&\frac{1}{(2\pi\Delta t)^{r/2}\sqrt{\det a(x)}}\times\\
&\exp\!\Bigl\{
-\tfrac{1}{2\Delta t}\|y-x\|_{a^{-1}(x)}^2\\
&-(y-x)\!\cdot a^{-1}(x)f(x)\\
&-\Delta t\bigl(\nabla\!\cdot f(x)+\tfrac12\|f(x)\|_{a^{-1}(x)}^2\bigr)
\Bigr\}
\end{aligned}
\) &
\(\displaystyle
\begin{aligned}
&\frac{1}{(2\pi\Delta t)^{r/2}\sqrt{\det_{+}a(x)}}\times\\
&\exp\!\Bigl\{
-\tfrac{1}{2\Delta t}\|y-x\|_{a^\dagger(x)}^2\\
&-(y-x)\!\cdot a^\dagger(x) f(x)\\
&-\Delta t\bigl(\nabla\!\cdot f(x)+\tfrac12\|f(x)\|_{a^\dagger(x)}^2\bigr)
\Bigr\}
\end{aligned}
\)
\\[2ex] \hline
\textbf{Online update} &
\(\displaystyle
\exp\!\Bigl\{
h^\top b^{-1}\,\Delta y_k
-\tfrac12\,\Delta t\,h^\top b^{-1}h
\Bigr\}\,\rho_k(x,\tau_k)
\) &
\(\displaystyle
\exp\!\Bigl\{
h^\top b^\dagger\,\Delta y_k
-\tfrac12\,\Delta t\,h^\top b^\dagger h
\Bigr\}\,\rho_k(x,\tau_k)
\)
\\
\bottomrule
\end{tabular}}
\end{table*}

\section{Numerical Divergence Calculation}
\label{sec:complex-step-divergence}

Let \(f:\mathbb{R}^r\to\mathbb{R}^r\) be a sufficiently smooth vector field,
\[
f(x)=(f_1(x),\dots,f_r(x))^\top,\quad x\in\mathbb{R}^r,
\]
and consider its divergence
\[
\nabla\cdot f(x)
=\sum_{k=1}^r \frac{\partial f_k(x)}{\partial x_k},
\]
which is required in the kernel construction \Cref{eq:kernel}.  Although closed‐form expressions for \(\nabla\cdot f\) can be used when available, deriving and implementing them can be cumbersome or even impossible for many complex applications. In such cases, one typically resorts to numerical approximation or automatic differentiation \cite{Baydin_2015}. Here we describe the complex‐step approximation \cite{martins_2003} for the \(k\)th partial derivative:
\[
\frac{\partial f_k}{\partial x_k}(x)
=\frac{\Im\bigl[f_k\bigl(x+i\,h\,e_k\bigr)\bigr]}{h}
+O(h^2),
\]
where
\begin{itemize}
  \item \(h>0\) is a small real step (e.g.\ \(h=10^{-20}\)),
  \item \(e_k\) is the \(k\)th standard basis vector in \(\mathbb{R}^r\),
  \item \(\Im(\cdot)\) denotes the imaginary part.
\end{itemize}
By summing over all components, the divergence can be approximated in only \(r\) evaluations of \(f\):
\[
\boxed{\nabla\cdot\!f(x)
\;\approx\;
\sum_{k=1}^r
\frac{\Im\bigl[f_k(x + i\,h\,e_k)\bigr]}{h}.}
\]
Compared with the classical central‐difference formula
\[
 \frac{\partial f_k}{\partial x_k}(x) \approx \frac{f_k(x+he_k)-f_k(x-he_k)}{2h},
\]
which requires $2r$ evaluations of $f$ for computing $\nabla \cdot f(x)$ and suffers from subtractive cancellation in the numerator when \(h\) is small (because it subtracts two nearly identical real values), the complex‐step approximation
avoids any subtraction of close real quantities and requires only $r$ evaluations of $f$. Consequently, one may take \(h\) extremely small (e.g.\ \(h\approx10^{-20}\)) and still recover derivatives to nearly machine‐precision accuracy.

\section{Log-Domain Computation}
\label{sec:log-domain-computation}

In both offline kernel evaluation and online filtering updates, numerous exponentiation and multiplication of tiny probability weights can cause numerical underflow or overflow. To enhance numerical stability, we perform all computations in the logarithmic domain. Let the sample set be \(\{x_j\}_{j=1}^n\), the discrete solution-mapping operator \(F_{\Delta t}\in\mathbb{R}^{n\times n}\), and the initial weights \(\sigma_0(x_j), j=1,\ldots,n\). The core formulas for log-domain computation are summarized below:

\paragraph{Prediction Step (Offline)} In the original domain,
\[
\sigma_{\mathrm{pred}}(i)
=
\sum_{j=1}^{n}F_{\Delta t}(j,i)\,\sigma_0(x_j).
\]
Hence, in the log domain,
\[
\log\sigma_{\mathrm{pred}}(i)
=
\operatorname{LSE}_{j}\bigl\{\log F_{\Delta t}(j,i)+\log\sigma_0(x_j)\bigr\},
\]
where the log‐sum‐exp operator is defined by
\[
\operatorname{LSE}_{j}\{a_j\}
=
\max_{j}a_j
\;+\;
\log\Bigl(\sum_{j}\exp(a_j-\max_{j}a_j)\Bigr).
\]

\paragraph{Observation Update (Online)}

Given the observation increment \(\Delta y\in \mathbb{R}^m\) and observation function \(h:\mathbb{R}^r \to \mathbb{R}^m\), the likelihood for sample \(x_i\) is
\[
L(i)
= \exp\Bigl(h(x_i)^{\top}\Delta y \;-\;\tfrac12\|h(x_i)\|^2\Delta t\Bigr),
\]
so the log‐weight update is simply
\[
\log\sigma_{\mathrm{upd}}(i)
=
\log\sigma_{\mathrm{pred}}(i)
+
h(x_i)^{\top}\Delta y - \frac{1}{2}\|h(x_i)\|^2\Delta t.
\]

\paragraph{Normalization}

To normalize weights without leaving the log domain, compute
\[
\log \tilde\sigma_{\mathrm{upd}}(i)
=
\log\sigma_{\mathrm{upd}}(i)
-
\operatorname{LSE}_{k}\{\log\sigma_{\mathrm{upd}}(k)\},
\]
so that
\[
\tilde\sigma_{\mathrm{upd}}(i)
=
\exp\bigl(\log \tilde\sigma_{\mathrm{upd}}(i)\bigr)
=
\frac{\sigma_{\mathrm{upd}}(i)}
     {\sum_{k=1}^n\sigma_{\mathrm{upd}}(k)}.
\]

Overall, carrying out all exponential and multiplicative operations in the log domain ensures robust numerical performance and avoids underflow/overflow in both offline kernel evaluation and online filtering updates.

\section{CPU/GPU Parallel Acceleration}
\label{sec:cpu-gpu-parallel-acceleration}

When the number of sampling points \(n\) becomes large, both the offline construction of the solution‐mapping matrix and the online update steps become computational bottlenecks.  We suggest accelerating these two stages as follows:

\paragraph{Offline: Multi‐threaded CPU assembly}
Building the \(n\times n\) operator matrix \(F_{\Delta t}\) entails \(\mathcal O(n^2)\) kernel evaluations.  To reduce wall‐clock time, we distribute the row‐wise kernel computations across multiple CPU threads or workers via parallel computing (e.g.\ MATLAB's \texttt{parfor} or an OpenMP implementation), yielding significant speedups for large \(n\) and high dimension $r$.

\paragraph{Online: GPU‐accelerated matrix-vector multiplies}
In the online update, the dominant cost is the repeated multiplication of the precomputed solution‐operator matrix by the weight vector.  We offload these intensive matrix-vector products to the GPU, which dramatically reduces per‐step runtime for large \(n\) and high dimension $m$.

Implementing the improved \yauyau filter with multi-threaded CPU matrix construction offline and GPU-accelerated updates online significantly boosts performance, enabling near‐real-time processing of large sample sizes and high-dimensional problems over long time horizons.

\section{Numerical Experiments}
\label{sec:numerical-experiments}

In this section, we validate our improved \yauyau filtering algorithm on three complementary testbeds:  
\begin{enumerate}
  \item[A.] \emph{Large‐scale scalability} (up to $r=1000$);  
  \item[B.] \emph{Small‐scale accuracy} with comparison to Extended Kalman Filter (EKF), Unscented Kalman Filter (UKF) and Particle Filter (PF);
  \item[C.] \emph{Linear system benchmark} versus the Kalman-Bucy filter.   
\end{enumerate}
\paragraph{Test Environment}
All experiments run in MATLAB R2023b on a workstation with  
\begin{itemize}
  \item Intel i9-12900K CPU (16 cores @ 3.7 GHz), 128 GB RAM, one NVIDIA RTX 4090 GPU;  
  \item MATLAB Parallel Toolbox (16 workers) for offline matrix construction;  GPU computing for online updates.  
\end{itemize}
\paragraph{Parameters}
Each experimental configuration is specified by:
\begin{itemize}
  \item \textbf{Dimensions \& sampling:} state dimension $r$, observation dimension $m$, and number of samples $n$;
  \item \textbf{Temporal grid:} time horizon $T$, number of steps $K$, and step size $\Delta t = T/K$;
  \item \textbf{Spatial sampling:} low‐discrepancy sequence (Halton for $r<10$, Sobol for $r\geq 10$);
  \item \textbf{Resampling-restart:} half‐width $R$ of the hypercube $[-R,R]^r$ (serving as the global range without resampling and the local range with resampling-restart), local samples $n_R$ and restart interval $\iota$.
\end{itemize}
\paragraph{Error Metrics} To quantify filter accuracy, we compute two error metrics over $K$ time steps:  
\begin{itemize}
    \item \textbf{Root-Mean-Square Error (RMSE):}
    \[
      \mathrm{RMSE}
      = \sqrt{\frac{1}{K\,r}\sum_{k=1}^{K}\|\hat x_k - x_k\|^2}\,.
    \]
    \item \textbf{Mean Error (ME):} 
    \[\mathrm{ME}  = \frac{1}{K}\sum_{k=1}^{K}\sqrt{\frac{1}{r}\|\hat x_k - x_k\|^2}.\]
\end{itemize}
Here, dividing by the state dimension $r$ removes the dependence on the number of state variables and yields the average per-component error. 
\begin{remark}
    A close match between RMSE and ME implies that estimation errors are uniformly distributed across time steps, reflecting robust and consistent performance. By contrast, if RMSE significantly exceeds ME, few large‐error steps dominate the mean‐square metric. Therefore, a small gap between RMSE and ME denotes stability and homogeneous errors, while a large gap highlights instability caused by extreme outliers.  
\end{remark}

\subsection{Large‐Scale Scalability Test}
\label{subsec:large-scale-test}

In this subsection, we assess the performance of the improved \yauyau filter on a suite of high‐dimensional, highly nonlinear problems as the state dimension $r$ increases from 10 to 1000.  We will demonstrate that the average RMSE and ME grows at most linearly, while the total runtime scales nearly linearly  in the state dimension $r$. This excellent behavior confirms that the improved \yauyau filter effectively mitigates the curse of dimensionality in high dimensions.

\paragraph{Setup}
Consider the cubic sensor system 
\begin{equation}\label{eq:cubic-sensor-system}
\begin{cases}
\mathrm{d}x_t = f(x_t)\,\mathrm{d}t+ d v_t,\\
\mathrm{d}y_t = h(x_t)\,\mathrm{d}t+ d w_t,
\end{cases}
\end{equation}
where \(x_t,y_t\in\mathbb{R}^r\), \(x_0\sim\mathcal{N}(0,I_r)\), and \(v_t,w_t\) are independent standard Brownian motions.  The drift function combines linear and periodic terms as
\begin{equation}
    f(x)=\sin(x)\odot (Ax)+\sin(2x)\odot (A_1x)
\end{equation}
where \(\odot\) denotes elementwise multiplication, and
\begin{equation}
    \label{eq:A&A1}
  A = \left[\begin{smallmatrix}
    -0.5 & 0.1 &        &  \\
         & -0.5 & \ddots &  \\
         &      & \ddots & 0.1 \\
         &      &        & -0.5
  \end{smallmatrix}\right],
  A_1 = \left[\begin{smallmatrix}
    -0.3 & 0.3 &        &  \\
         & -0.3 & \ddots &  \\
         &      & \ddots & 0.3 \\
         &      &        & -0.3
  \end{smallmatrix}\right].
\end{equation}
The observation function is
\[
  h(x) = (x - 100)^3.
\]
Experiments are conducted in the high‐dimensional regime:
\begin{itemize}
  \item \textbf{Dimension:} \(r\in\{10,50,100,300,600,1000\}\) and $m=r$.
  \item \textbf{Temporal grid:} \(T=10\), \(K=1000\) steps (\(\Delta t=0.01\)).
  \item \textbf{Spatial sampling:} Sobol sampling with sample points \(n=100, 300, 500, 800, 1000, 2000\) for dimensions \(r=10, 50, 100, 300, 600, 1000\) respectively.
  \item \textbf{Resampling-restart:} $R=0.3$, $n_R = n$ and $\iota=2$.
  \item \textbf{Monte Carlo trials:} \(20\) runs.
\end{itemize}
Both the true state trajectories and the observation data are generated by discretizing the underlying SDEs using the Euler–Maruyama scheme with step size $\Delta t$.

\paragraph{Results}  

\begin{table}[htbp]
\centering
\caption{Computation time, RMSE, and ME of the improved Yau–Yau filter for large-scale cases under varying dimensions \(r\) and sample sizes \(n_{R}\).  Values are mean \(\pm\) standard deviation over 20 runs.}
\label{tab:case6_performance}
\begin{tabular}{ccccc}
\toprule
$r$ & $n_{R}$ & Time (s) & RMSE & ME \\
\midrule
10   & 100   & $0.271 \pm 0.023$   & $0.165 \pm 0.002$ & $0.160 \pm 0.002$ \\
50   & 300   & $1.006 \pm 0.026$   & $0.368 \pm 0.001$ & $0.365 \pm 0.001$ \\
100  & 500   & $4.824 \pm 0.095$   & $0.523 \pm 0.005$ & $0.520 \pm 0.005$ \\
300  & 800   & $16.833 \pm 0.374$  & $0.964 \pm 0.034$ & $0.955 \pm 0.033$ \\
600  & 1000  & $37.160 \pm 0.065$  & $1.376 \pm 0.049$ & $1.356 \pm 0.046$ \\
1000 & 2000  & $150.540 \pm 6.431$ & $1.609 \pm 0.013$ & $1.582 \pm 0.012$ \\
\bottomrule
\end{tabular}
\end{table}

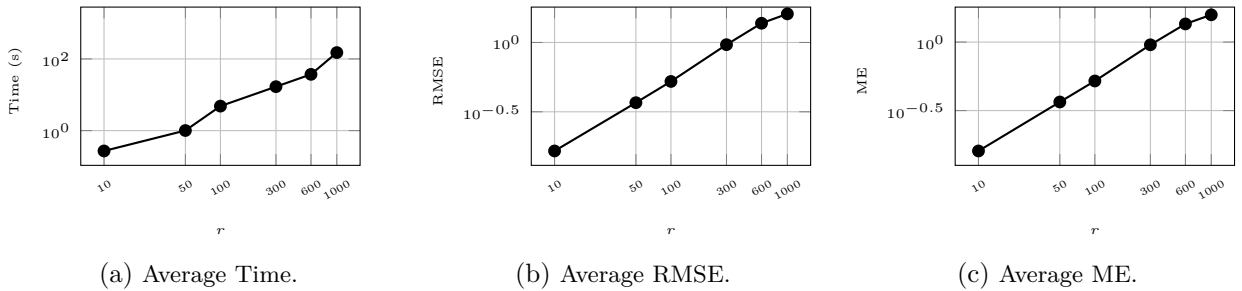
\begin{figure*}[ht]
\centering
\begin{subfigure}[t]{0.32\textwidth}
\begin{tikzpicture}
\begin{axis}[
    xmode=log,
    ymode=log,
    log basis x=10,
    log basis y=10,
    xlabel={\(r\)},
    ylabel={Time (s)},
    xtick={10,50,100,300,600,1000},
    xticklabels={10,50,100,300,600,1000},
    grid=both,
    width=\linewidth,
    height=0.7\linewidth,
    ymin=0, ymax=2800,
    xticklabel style={font=\fontsize{4pt}{4pt}\selectfont, rotate=30},
    yticklabel style={
      font=\fontsize{4pt}{4pt}\selectfont
    },
    xlabel style={
      font=\fontsize{5pt}{5pt}\selectfont,
      yshift=-2pt
    },
    ylabel style={
      font=\fontsize{5pt}{5pt}\selectfont,
      xshift=2pt
    },
]
\addplot[mark=*, thick] coordinates {
    (10,0.27051)
    (50,1.0062)
    (100,4.8244)
    (300,16.8330)
    (600,37.1600)
    (1000,150.5400)
};
\end{axis}
\end{tikzpicture}
\caption{{\footnotesize Average Time.}}
\label{fig:case6_time_new}
\end{subfigure}
\hfill
\begin{subfigure}[t]{0.32\textwidth}
\begin{tikzpicture}
\begin{axis}[
    xmode=log,
    ymode=log,
    log basis x=10,
    log basis y=10,
    xlabel={\(r\)},
    ylabel={RMSE},
    xtick={10,50,100,300,600,1000},
    xticklabels={10,50,100,300,600,1000},
    grid=both,
    width=\linewidth,
    height=0.7\linewidth,
    ymin=0, ymax=1.8,
    xticklabel style={
     font=\fontsize{4pt}{4pt}\selectfont, rotate=30
    },
    yticklabel style={
      font=\fontsize{4pt}{4pt}\selectfont
    },
    xlabel style={
      font=\fontsize{5pt}{5pt}\selectfont,
      yshift=-2pt
    },
    ylabel style={
      font=\fontsize{5pt}{5pt}\selectfont,
      xshift=2pt
    },
]
\addplot[mark=*, thick] coordinates {
    (10,0.16497)
    (50,0.36787)
    (100,0.52337)
    (300,0.96410)
    (600,1.37600)
    (1000,1.60850)
};
\end{axis}
\end{tikzpicture}
\caption{{\footnotesize Average RMSE.}}
\label{fig:case6_rmse_new}
\end{subfigure}
\hfill
\begin{subfigure}[t]{0.32\textwidth}
\begin{tikzpicture}
\begin{axis}[
    xmode=log,
    ymode=log,
    log basis x=10,
    log basis y=10,
    xlabel={\(r\)},
    ylabel={ME},
    xtick={10,50,100,300,600,1000},
    xticklabels={10,50,100,300,600,1000},
    grid=both,
    width=\linewidth,
    height=0.7\linewidth,
    ymin=0, ymax=1.8,
    xticklabel style={
     font=\fontsize{4pt}{4pt}\selectfont, rotate=30
    },
    yticklabel style={
      font=\fontsize{4pt}{4pt}\selectfont
    },
    xlabel style={
      font=\fontsize{5pt}{5pt}\selectfont,
      yshift=-2pt
    },
    ylabel style={
      font=\fontsize{5pt}{5pt}\selectfont,
      xshift=2pt
    },
]
\addplot[mark=*, thick] coordinates {
    (10,0.16049)
    (50,0.36541)
    (100,0.52004)
    (300,0.95541)
    (600,1.35620)
    (1000,1.58150)
};
\end{axis}
\end{tikzpicture}
\caption{{\footnotesize Average ME.}}
\label{fig:case6_me}
\end{subfigure}
\caption{Performance metrics for large-scale nonlinear filtering over 1000 time steps as state dimension \(r\) varies, plotted on log-log axes: (a) average total computation time; (b) average RMSE; (c) average ME.}
\label{fig:case6_r_vs_metrics}
\end{figure*}

From Table~\ref{tab:case6_performance}, we observe that even with a modest local sample size (\(n_R=100\)–\(2000\)), the improved \yauyau filter maintains high accuracy: for \(r\le100\), both RMSE and ME remain below 0.6, and at \(r=1000\) they stay under 1.6, demonstrating reliable state estimation with limited sampling overhead.

As shown in Fig.~\ref{fig:case6_r_vs_metrics}, the total computation time scales nearly linearly in $r$ (approximately \(\mathcal{O}(r^{1.2})\)), while RMSE and ME increase only sub‐linearly with state dimension. Such polynomial growth is vastly smaller than the exponential blow-up one would expect in other classical methods, demonstrating effective mitigation of the curse of dimensionality.

\begin{figure*}[ht]
  \centering
  \begin{subfigure}[t]{0.45\linewidth}
    \includegraphics[width=\linewidth]{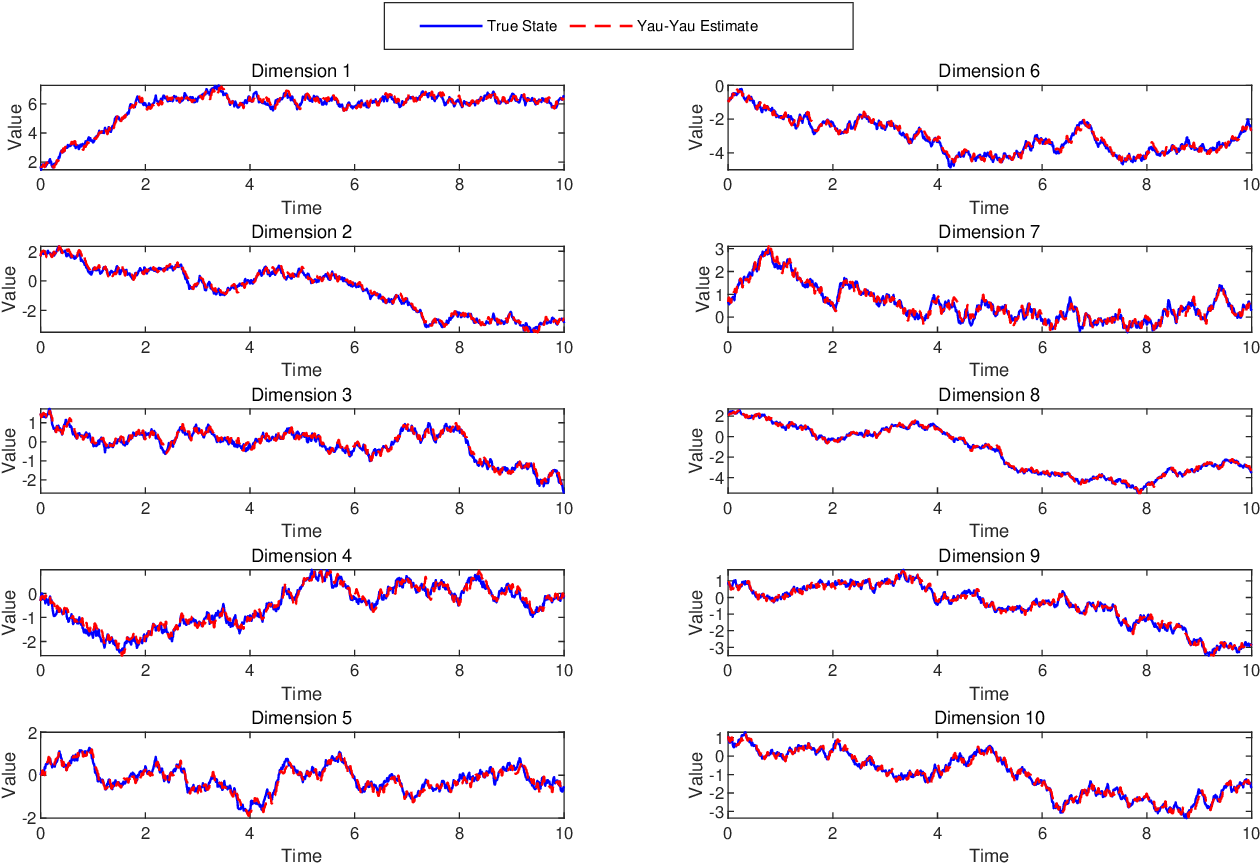}
    \caption{$r=10,\ n_{R}=100$}
    \label{fig:case6_r10}
  \end{subfigure}
  \hfill
  \begin{subfigure}[t]{0.45\linewidth}
    \includegraphics[width=\linewidth]{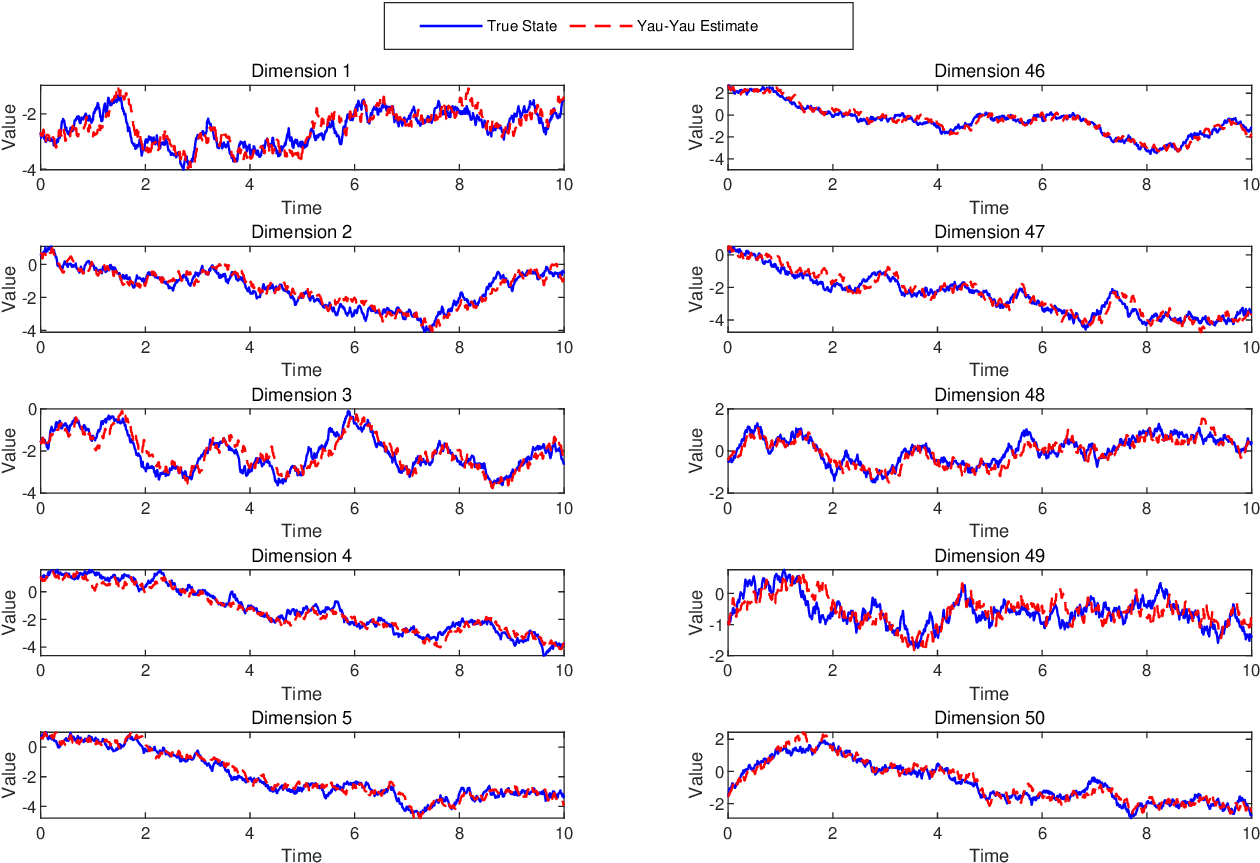}
    \caption{$r=50,\ n_{R}=300$}
    \label{fig:case6_r50}
  \end{subfigure}
  \hfill
  \begin{subfigure}[t]{0.45\linewidth}
    \includegraphics[width=\linewidth]{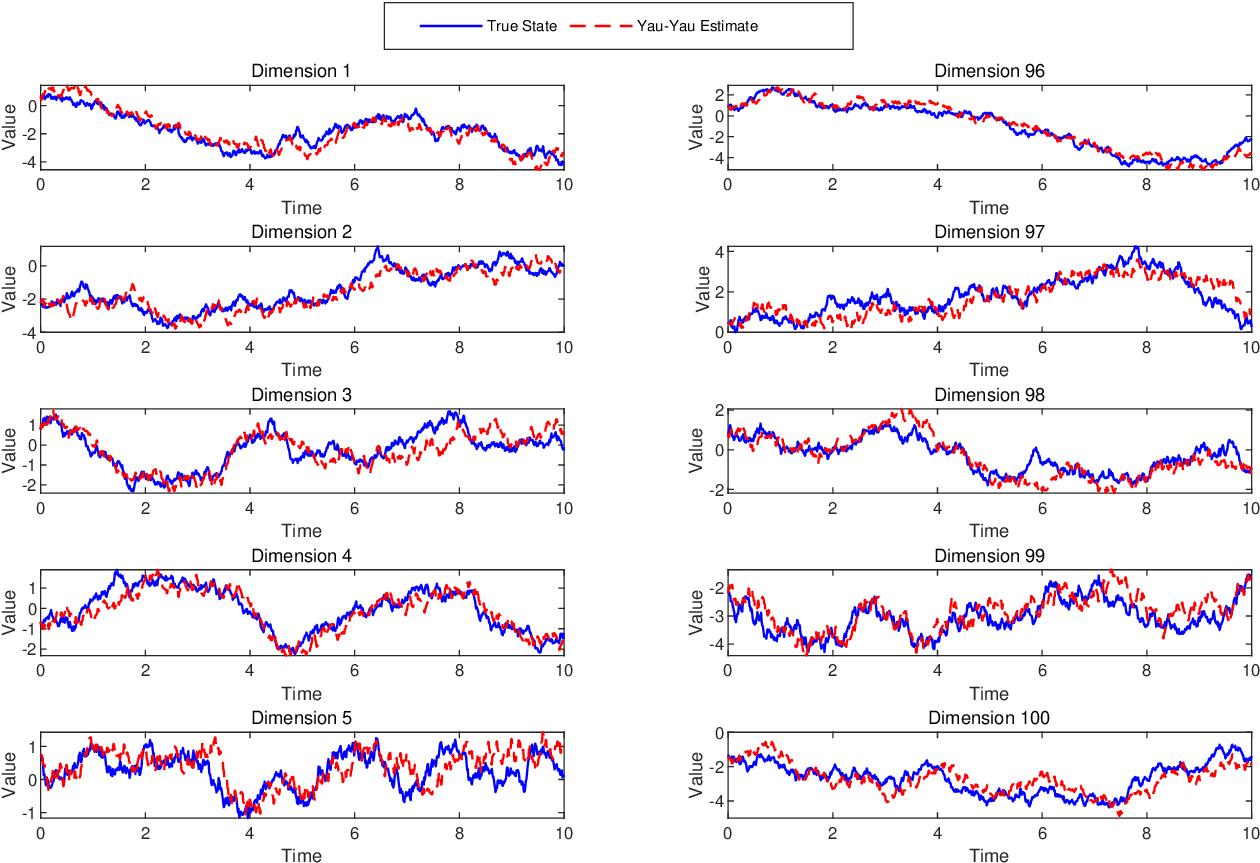}
    \caption{$r=100,\ n_{R}=500$}
    \label{fig:case6_r100}
  \end{subfigure}
    \caption{Improved Yau–Yau filtering results for the large-scale nonlinear filtering over 1000 time steps: 
  (a) $r=10$ with $n_R=100$; 
  (b) $r=50$ with $n_R=300$; 
  (c) $r=100$ with $n_R=500$. 
  In each subfigure, only the first five and last five state variables are shown, with the solid blue line indicating the true state trajectory and the dashed red line the improved Yau–Yau estimate.}

  \label{fig:case6_comparison}
\end{figure*}

Fig.~\ref{fig:case6_comparison} further illustrates that, over 1000 time steps with \(\Delta t=0.01\), the improved \yauyau filter accurately tracks the true state trajectories across multiple high-dimensional scenarios.

Overall, these results demonstrate that the improved \yauyau algorithm is both scalable and accurate in high-dimensional settings, effectively overcoming the curse of dimensionality. 

\subsection{Small‐Scale Accuracy Comparison}
\label{subsec:small-scale}

In this subsection, we evaluate four nonlinear filtering methods—EKF, UKF, PF, and the Improved \yauyau Filter—on two challenging small‐scale test problems.

The first test uses a highly nonlinear observation function \(h(x)=1000x^3\).  Here EKF and UKF suffer catastrophic degradation due to the failure of local linearization around zero, while PF and the Improved \yauyau filter both succeed in tracking the true state.

The second test is the “double‐well” potential with square observation \(h(x)=x^2\).  Initialized at the saddle point \(x_0=0\), EKF and UKF produce zero Kalman gain and remain stuck at the initial estimate, PF with limited particles under-samples one well and exhibits high variance, whereas the improved \yauyau filter remains robust and delivers the best accuracy and stability across multiple trials.

\subsubsection{Small-Scale Highly Nonlinear Test}
\paragraph{Setup}  
We first consider the one-dimensional model
\[
\mathrm{d}x_t = -\,x_t\,\mathrm{d}t + \mathrm{d}v_t,
\qquad
\mathrm{d}y_t = 1000 x_t^3\,\mathrm{d}t + \mathrm{d}w_t.
\]
The large coefficient \(1000\) makes the observation function extremely sensitive—even small deviations in \(x\) produce dramatic changes in \(h(x) \;=\; 1000\,x^3\).  Moreover, the non‐shifted cubic observation function
\(h(x)\) cannot be effectively linearized around zero, causing EKF and UKF to fail catastrophically in some neighborhood of \(0\).

We simulate over $T=10$ with the time step $\Delta t=0.01$ (total $1000$ steps) and the initial point $x_0 = 0.5$. For both particle‐based methods (PF and improved \yauyau), we set $n=200$ samples; additionally, the improved \yauyau filter employs a local resampling radius $R=4.5$ with $n_R=200$ sample points, and performs resampling every $\iota=16$ steps.

\paragraph{Result} 
\begin{table}[ht]
\centering
\caption{Performance comparison on the small-scale test.}
\label{tab:filter_performance_sp1}
\begin{tabular}{lccc}
\hline
Method & RMSE & ME & Time (s) \\
\hline
EKF & 0.5967 & 0.1679 & 0.0016 \\
UKF & 0.2078 & 0.1345 & 0.0033 \\
PF & \textbf{0.0626} & \textbf{0.0342} & 0.0454 \\
Yau–Yau & 0.0645 & 0.0402 & 0.2273 \\
\hline
\end{tabular}
\end{table}

\begin{figure*}[ht!]
  \centering
  \begin{subfigure}[t]{0.48\textwidth}
    \centering
    \includegraphics[width=\textwidth]{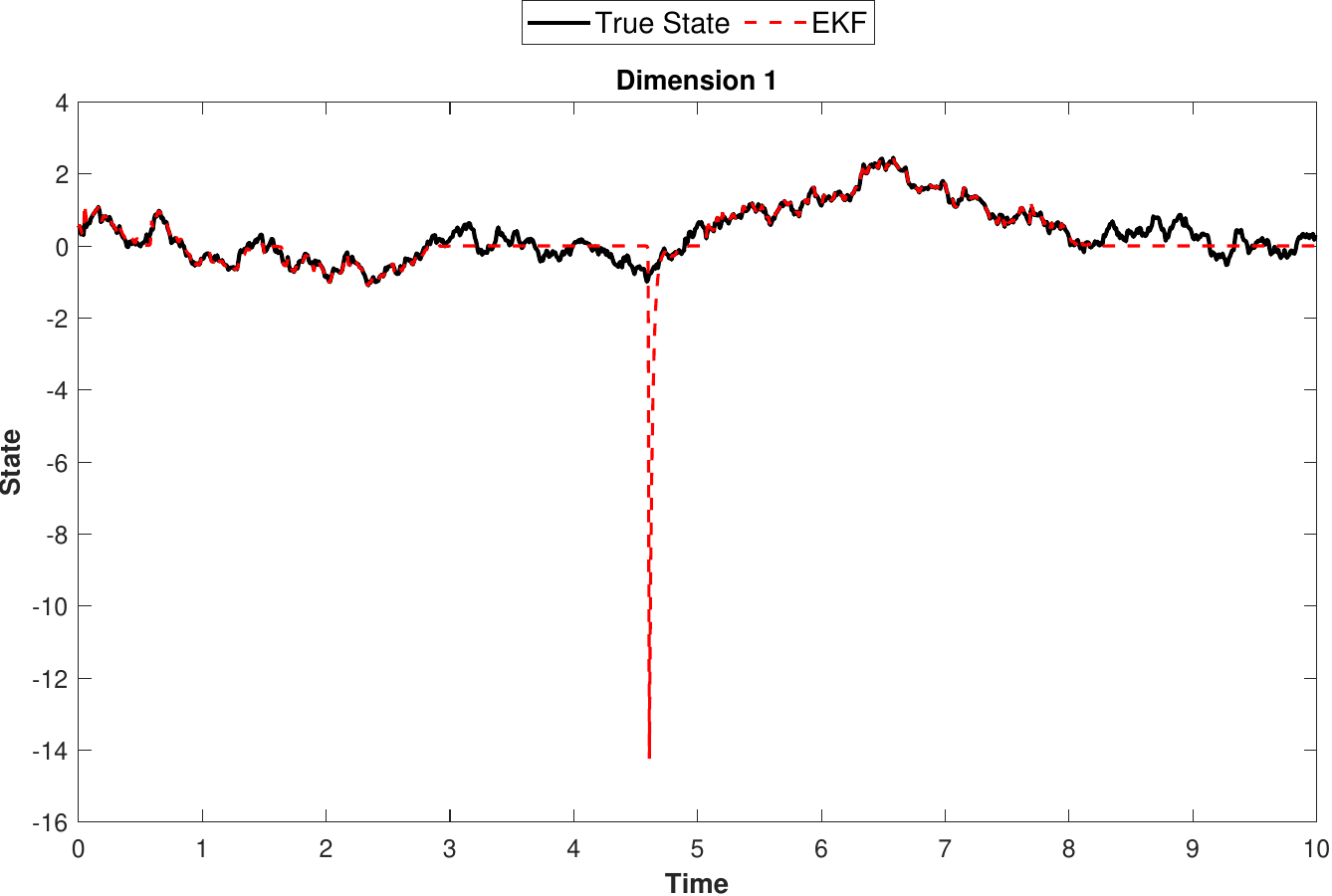}
    \caption{EKF}
  \end{subfigure}%
  \hfill
  \begin{subfigure}[t]{0.48\textwidth}
    \centering
    \includegraphics[width=\textwidth]{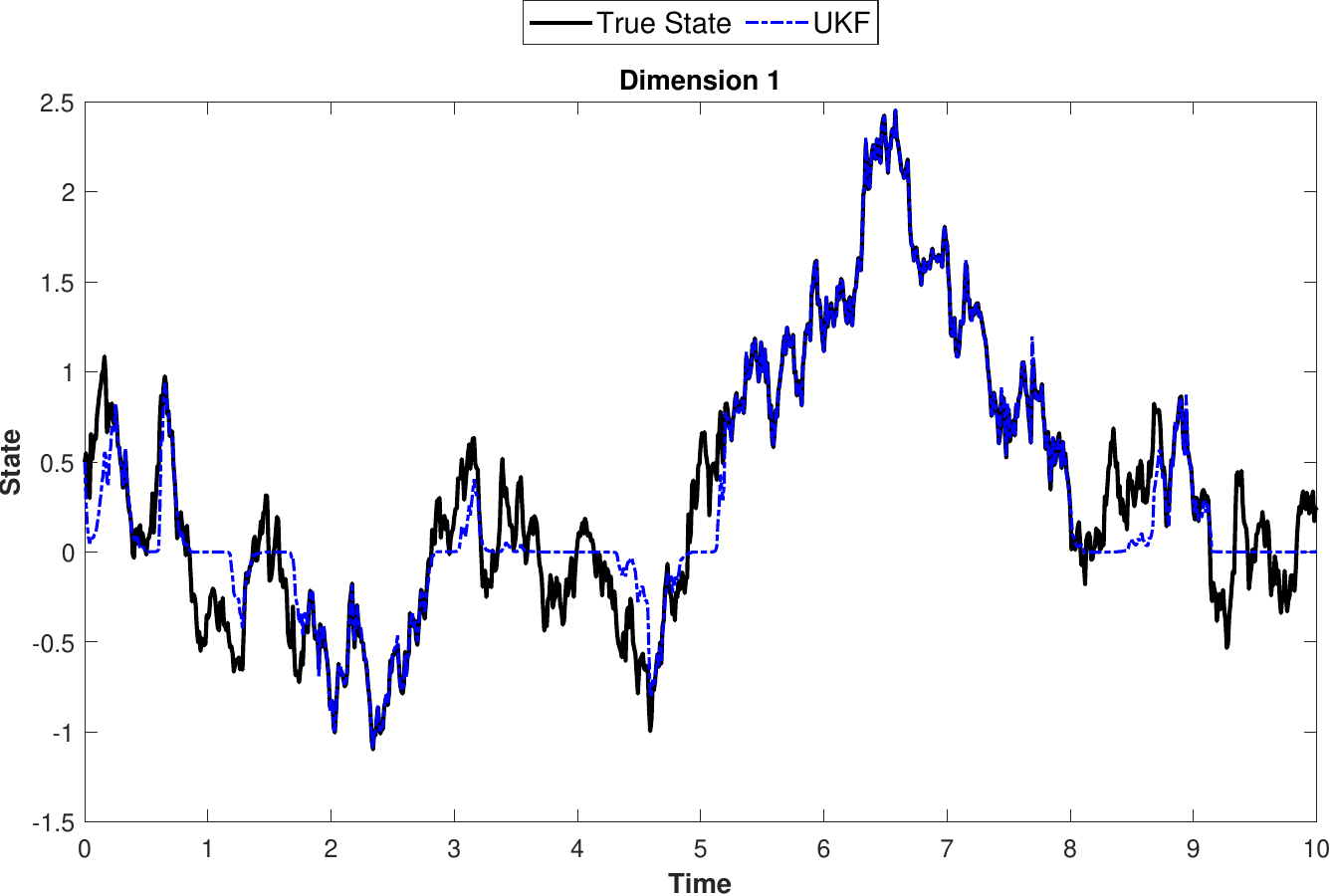}
    \caption{UKF}
  \end{subfigure}%
  \hfill
  \begin{subfigure}[t]{0.48\textwidth}
    \centering
    \includegraphics[width=\textwidth]{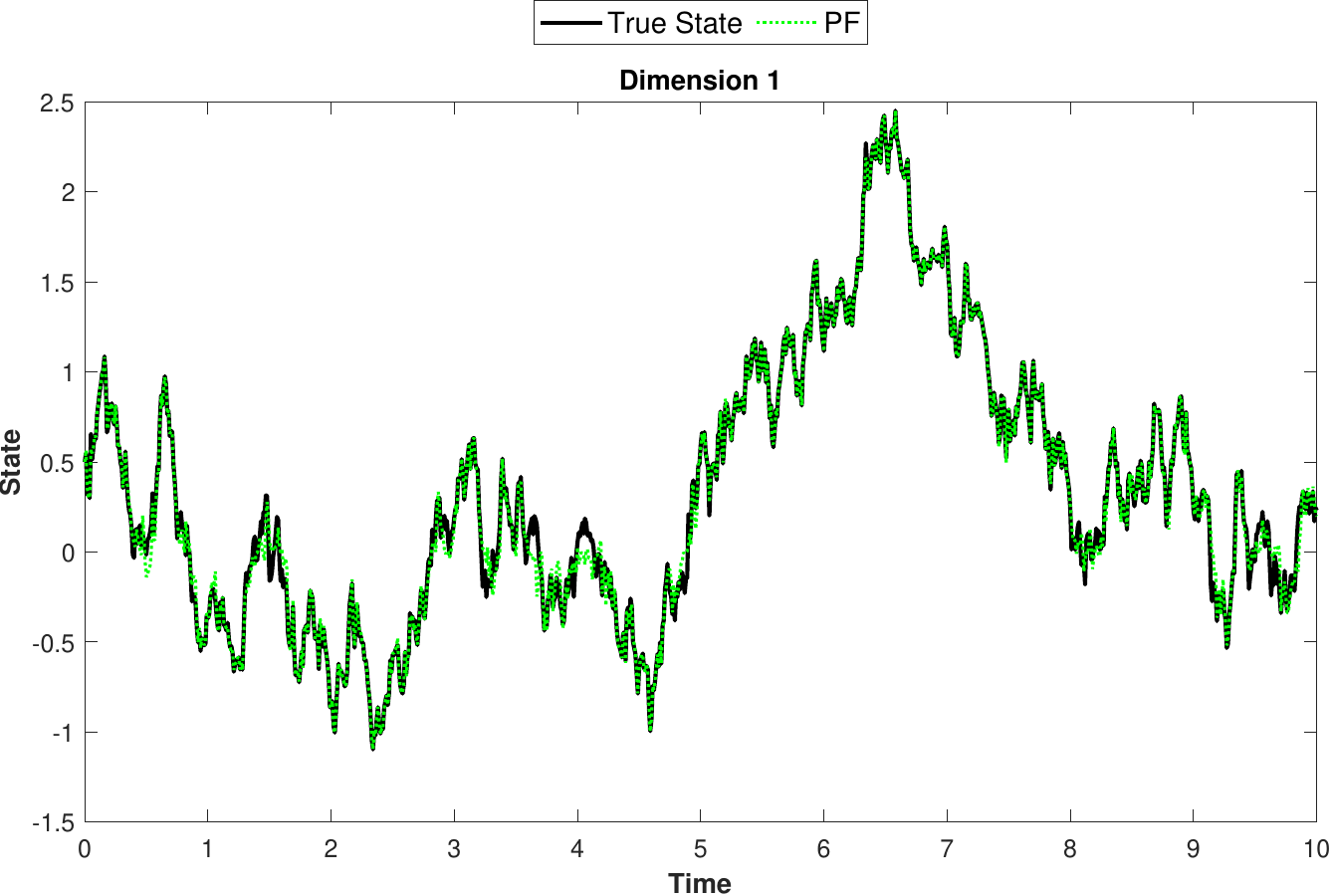}
    \caption{PF}
  \end{subfigure}%
  \hfill
  \begin{subfigure}[t]{0.48
  \textwidth}
    \centering
    \includegraphics[width=\textwidth]{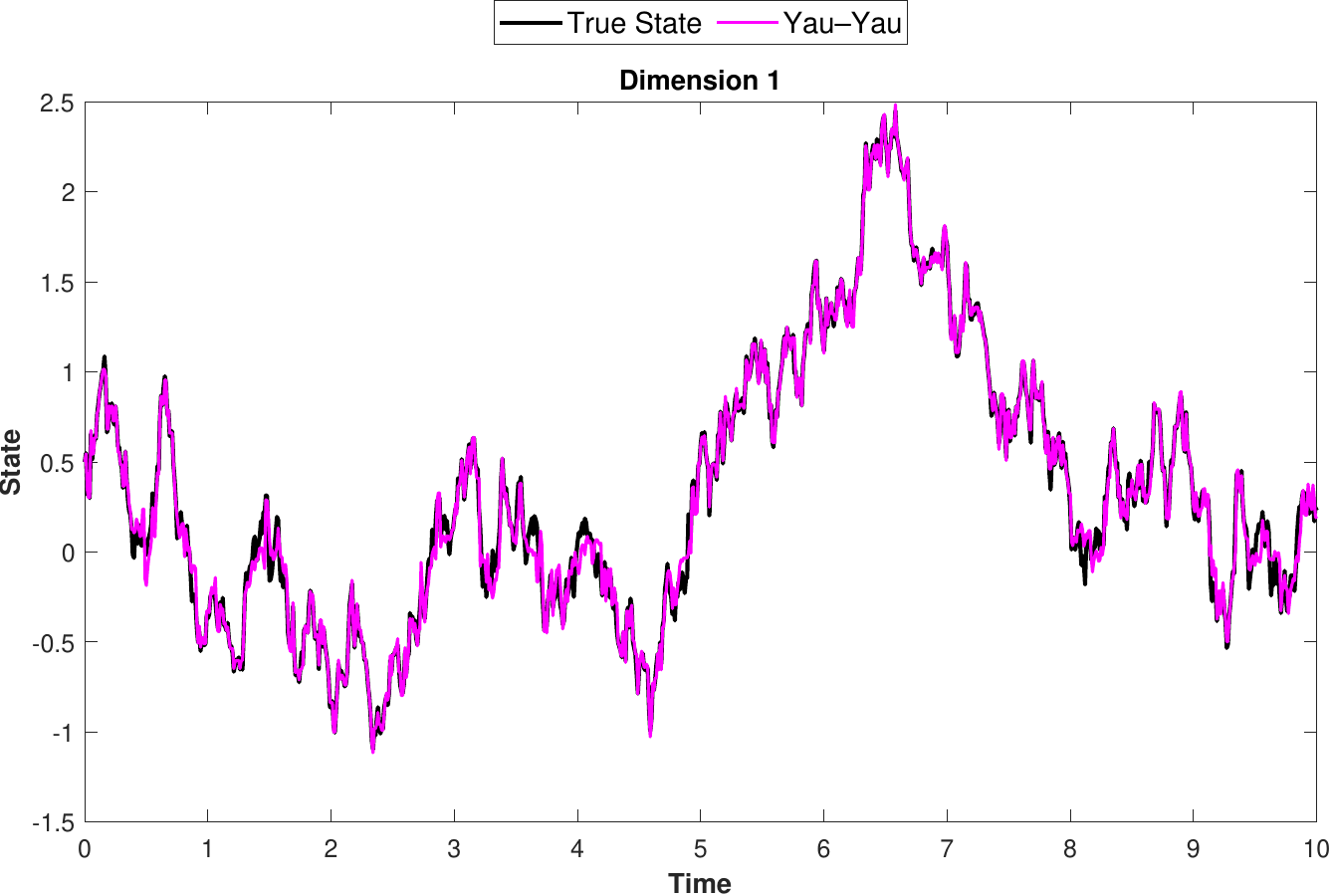}
    \caption{Yau–Yau}
  \end{subfigure}
  \caption{Performance comparison of four filtering algorithms on the 1-D cubic sensor test problem.  
    Each subpanel shows the time‐series estimate of the single state variable: (a) EKF; (b) UKF; (c) PF; (d) Yau–Yau.  
    All simulations use the same time-discretization ($\Delta t=0.01$, $T=10$, total $1000$ steps) and standard Gaussian noise settings.}
  \label{fig:filter_comparison}
\end{figure*}

As shown in Table~\ref{tab:filter_performance_sp1} and Fig.~\ref{fig:filter_comparison}, the improved \yauyau algorithm, while slightly more expensive in computation time (0.2273 s), produces the second-best state estimates (RMSE 0.0645 and ME 0.0402) among the four compared methods. By contrast, the EKF and UKF suffer severe instability around zero, resulting in substantially larger errors. The PF, however, performs well, achieving slightly better accuracy comparable to the improved \yauyau algorithm; yet PF can suffer from high variance and erratic behavior in some strongly nonlinear double‐well regimes, while the improved \yauyau filter not only outperforms PF in accuracy but also maintains significantly better stability, as we will see in the next test case.

\subsubsection{Double‐Well Potential with Square Observation Test}
\paragraph{Setup}
Consider the one‐dimensional “double‐well” SDE
\[
  \mathrm{d}x_t = -4\,x_t\bigl(x_t^2 - 1\bigr)\,\mathrm{d}t + \sigma\,\mathrm{d}v_t,\qquad
  \mathrm{d}y_t = x_t^2\,\mathrm{d}t + \rho\,\mathrm{d}w_t,
\]
with $\sigma=0.5$, $\rho=0.2$, discretized by Euler–Maruyama over $T=5$ with $\Delta t=0.01$ (500 steps).  All filters are run for 20 independent trials.  For PF we use $300$ particles.  For the improved \yauyau filter we likewise use $n=300$ local samples drawn via Halton sequences, with fixed radius $R=10$ and no-resampling.

\begin{figure*}[ht!]
  \centering
  \begin{subfigure}[t]{0.3\textwidth}
    \includegraphics[width=\textwidth]{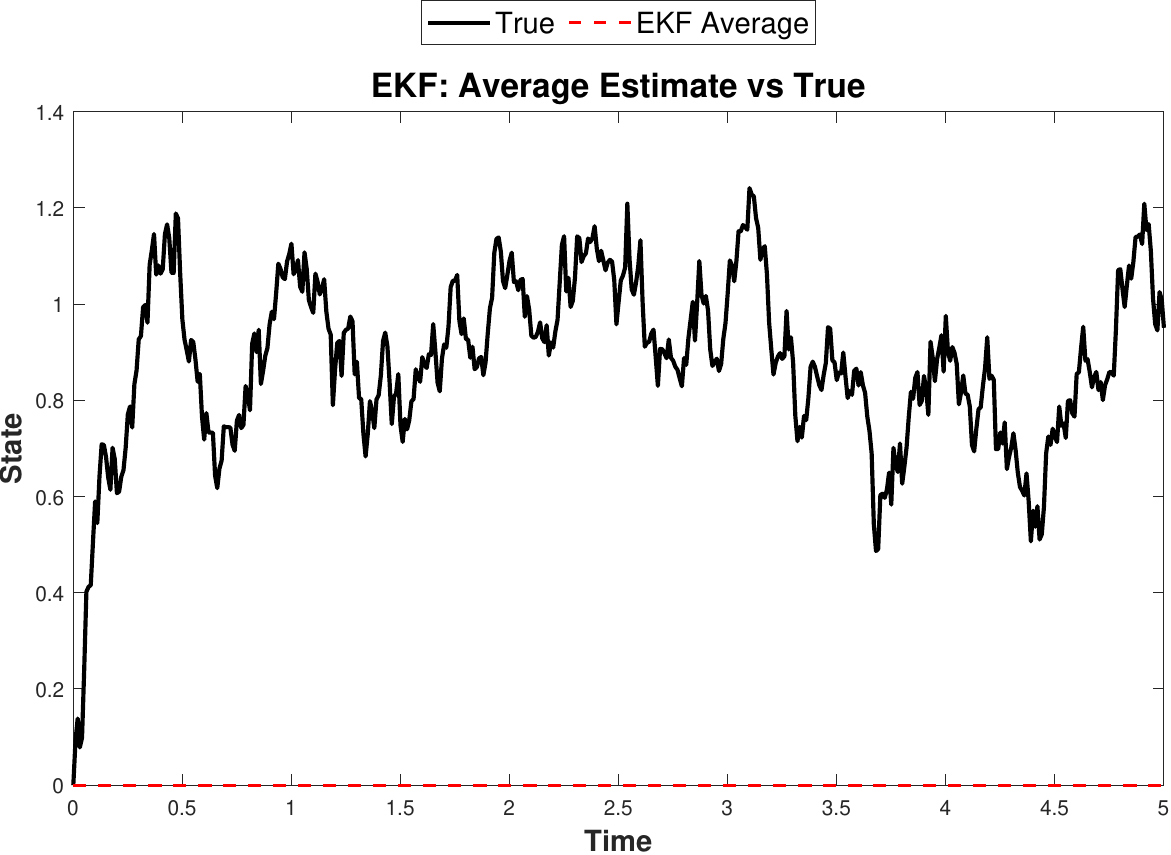}
    \caption{EKF}
  \end{subfigure}\quad
  \begin{subfigure}[t]{0.3\textwidth}
    \includegraphics[width=\textwidth]{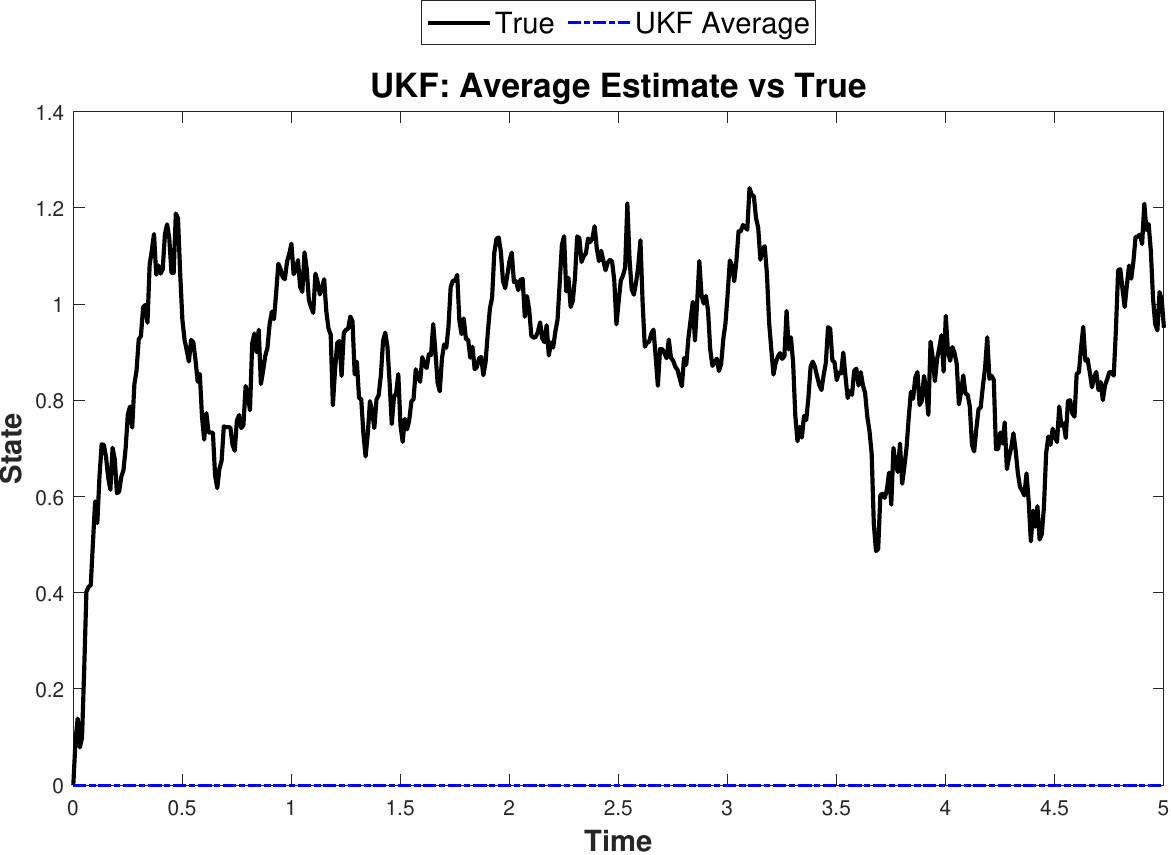}
    \caption{UKF}
  \end{subfigure}\quad
  \begin{subfigure}[t]{0.3\textwidth}
    \includegraphics[width=\textwidth]{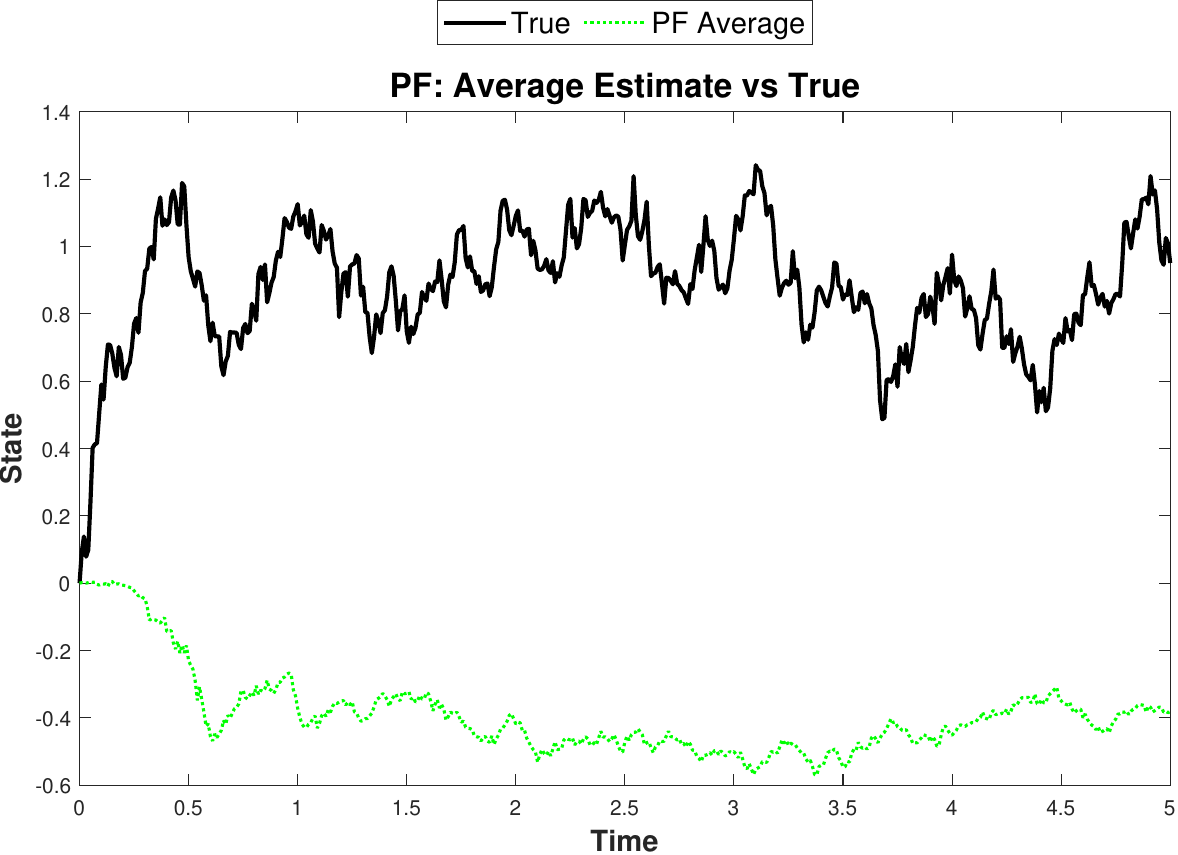}
    \caption{PF}
  \end{subfigure}\quad
  \begin{subfigure}[t]{0.3\textwidth}
    \includegraphics[width=\textwidth]{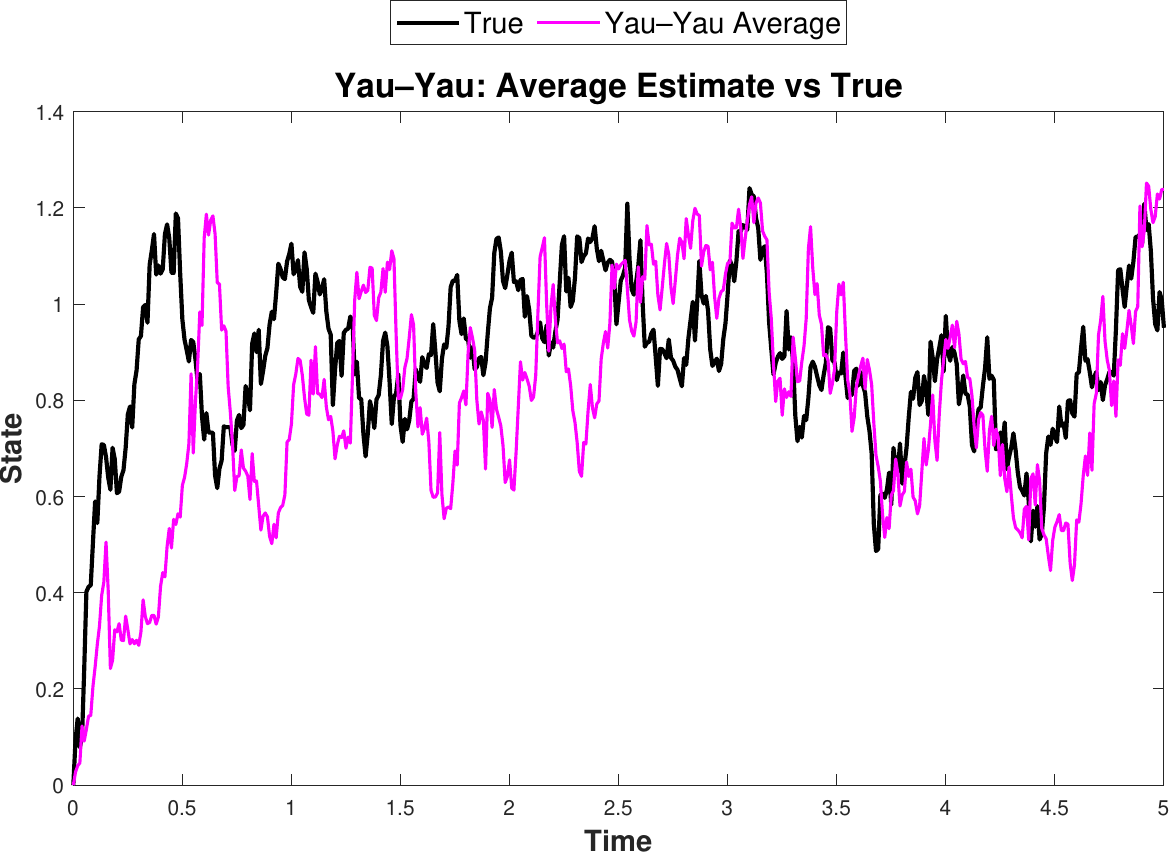}
    \caption{Yau–Yau}
  \end{subfigure}\quad
  \begin{subfigure}[t]{0.3\textwidth}
    \includegraphics[width=\textwidth]{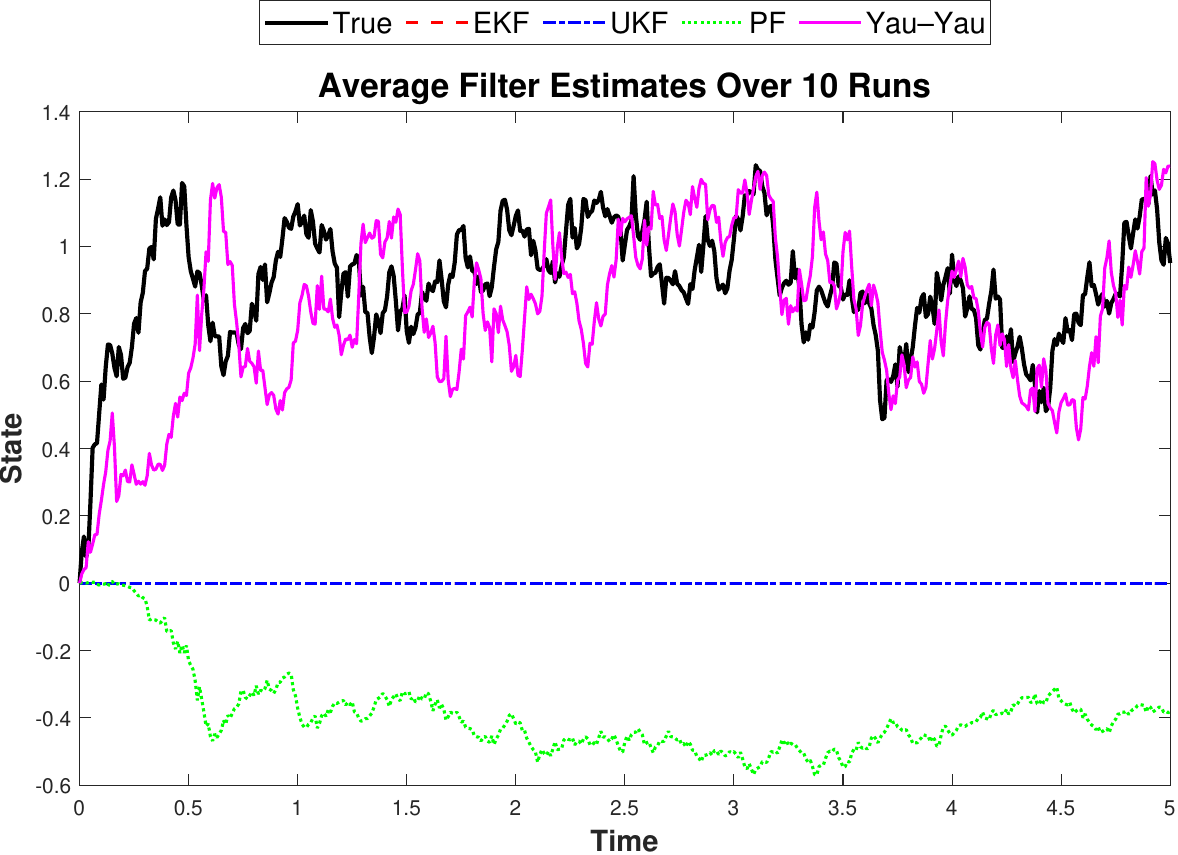}
    \caption{All methods}
  \end{subfigure}
  \caption{Average state estimates over 20 runs on the double‐well test:  
    (a) EKF and (b) UKF both diverge at the unstable equilibrium;  
    (c) PF collapses onto one well due to particle degeneracy;  
    (d) improved Yau–Yau accurately tracks both modes;  
    (e) overlay of all four against the true trajectory.}
  \label{fig:filter_comparison_double_well}
\end{figure*}

\paragraph{Results}
Table~\ref{tab:filter_double_well} and Fig.~\ref{fig:filter_comparison_double_well} summarize the performance. outcomes. Starting at the unstable equilibrium $x_0=0$, the state quickly falls into one of the wells at $x=\pm1$, yielding a strongly bimodal posterior under the square sensor $h(x)=x^2$. The EKF and UKF, relying on the linearization around the saddle point \(x=0\), yield zero Kalman gain and thus remain stuck at the initial estimate. The PF with only 300 particles severely under‐samples one mode, leading to collapse and high variance in its estimates.  In contrast, the improved \yauyau filter, using the same sample budget ($n=300$), remains both accurate and stable.  This clearly demonstrates the robustness of the improved \yauyau algorithm for strongly nonlinear, multimodal filtering problems.  

\begin{table}[ht]
  \centering
  \caption{Statistics over 20 runs on the double‐well test.}
  \label{tab:filter_double_well}
  \begin{tabular}{l c c c}
    \hline
    Method     & RMSE & ME    & Time (s) \\
    \hline
    EKF        & $0.9007\pm0.0000$ & $0.8826\pm0.0000$ & $0.0003\pm0.0003$ \\
    UKF        & $0.9007\pm0.0000$ & $0.8826\pm0.0000$ & $0.0014\pm0.0003$ \\
    PF         & $1.3603\pm0.4617$ & $1.2825\pm0.4990$ & $0.0141\pm0.0022$ \\
    Yau–Yau    & $\mathbf{0.2473\pm0.0000}$ & $\mathbf{0.1945\pm0.0000}$ & $0.1101\pm0.0053$ \\
    \hline
  \end{tabular}
\end{table}

\subsection{Linear Filtering Benchmark}
\label{subsec:linear-kalman}

We evaluate the improved \yauyau filter against the Kalman–Bucy filter on linear Gaussian filtering problems to demonstrate that our nonlinear solver remains competitive—and can even achieve superior—accuracy in linear settings.

\paragraph{Setup}  
Consider the continuous‐time linear model
\[
\mathrm{d}x = A\,x\,\mathrm{d}t + \mathrm{d}v_t,\qquad
\mathrm{d}y = 5\,x\,\mathrm{d}t + \mathrm{d}w_t,
\]
with \(A\) as defined in \eqref{eq:A&A1}, and \(v_t,w_t\) independent standard Wiener processes.  Simulations run over \(T=10\) with time step \(\Delta t=0.01\).  We use particle counts \(n=\{300,1000,1500\}\) for state dimensions \(r=\{1,2,3\}\), respectively, without resampling-restart. 

\begin{figure*}[htbp]
\centering
\includegraphics[width=0.32\linewidth]{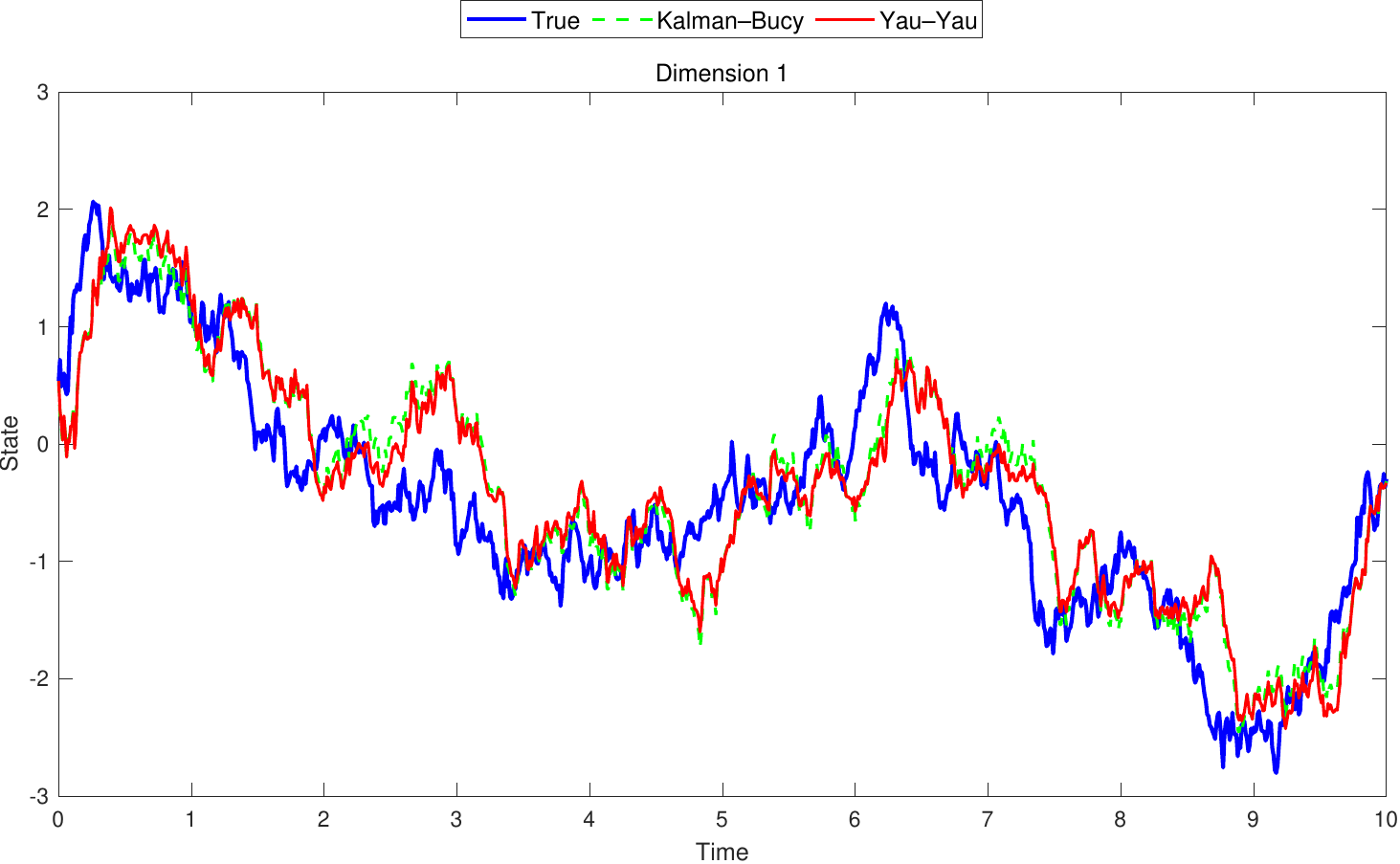}
\includegraphics[width=0.32\linewidth]{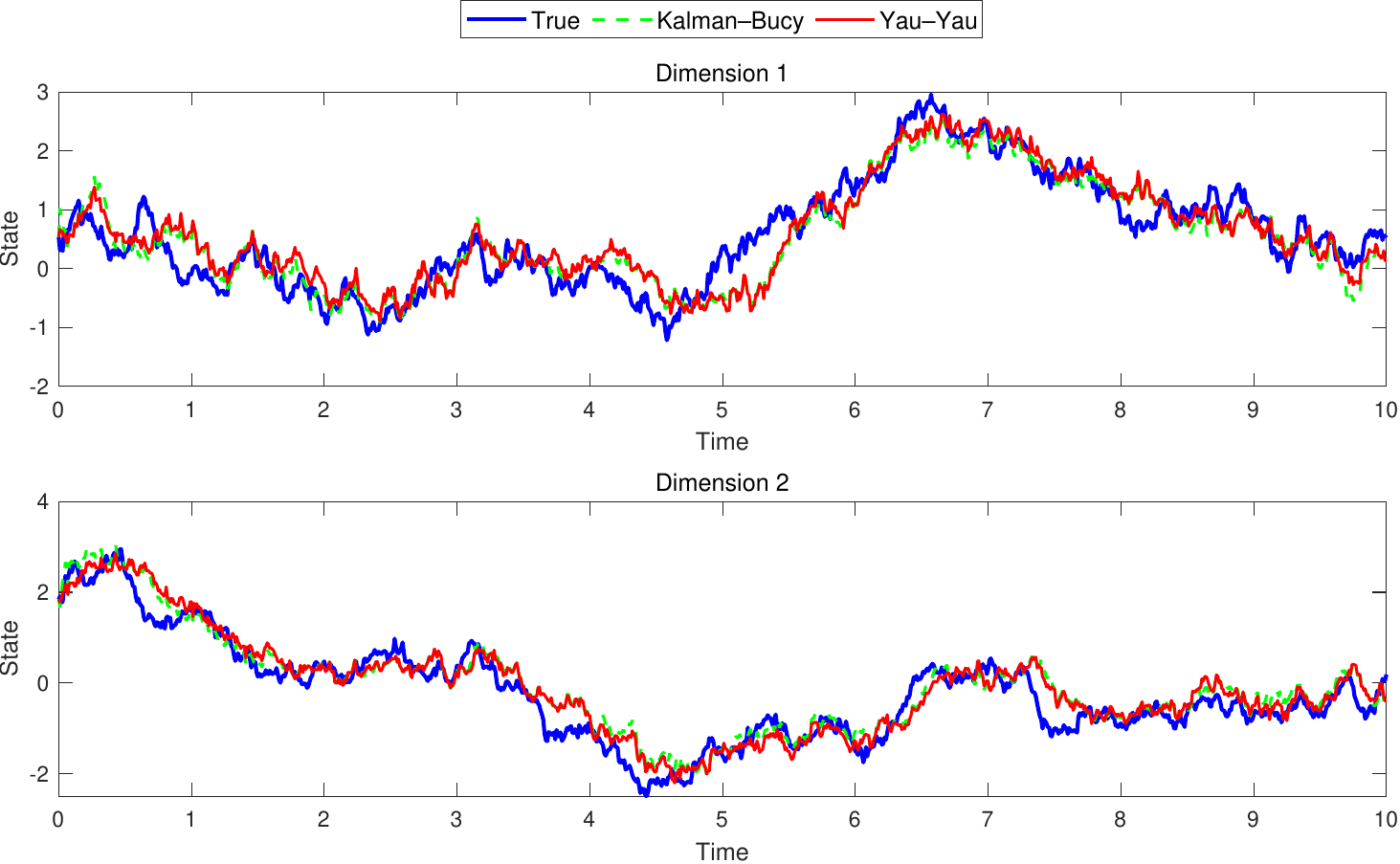}
\includegraphics[width=0.32\linewidth]{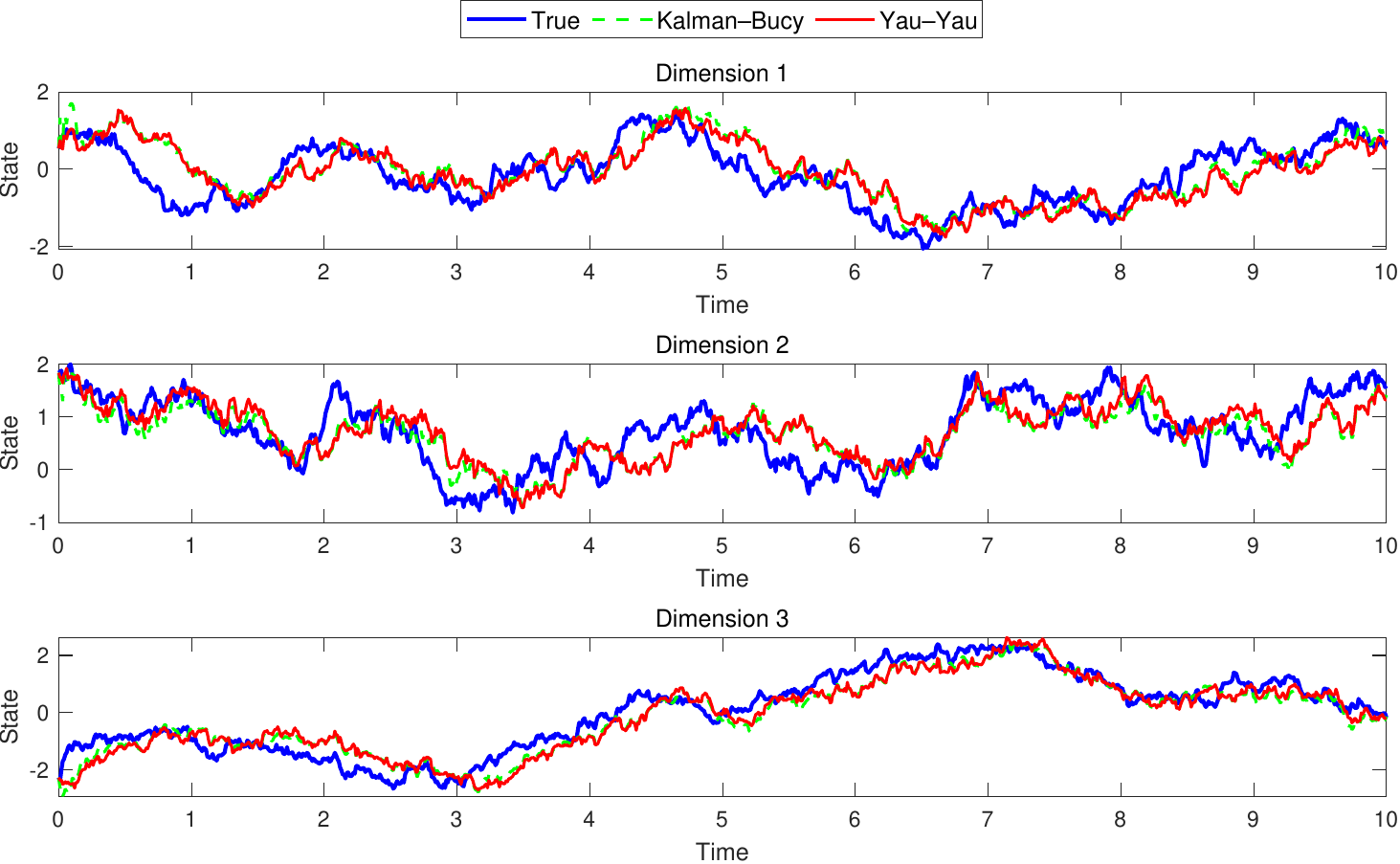}
\caption{State estimation trajectories for the improved Yau–Yau (dashed) and Kalman–Bucy (dotted) filters compared to the true state (solid) at \(r=1,2,3\) for linear filtering problems.}
\label{fig:linear_comparison}
\end{figure*}

\paragraph{Results}  
Table~\ref{tab:linear_performance} reports the RMSE, ME, and runtime for both filters. Surprisingly, the improved \yauyau filter matches or even outperforms the optimal Kalman–Bucy filter in estimation accuracy for \(r=1,2\), although Kalman–Bucy completes in under \(0.001\) seconds, \yauyau finishes in a modest \(0.25\)–\(2.64\) seconds. Fig.~\ref{fig:linear_comparison} shows state estimation trajectories: both filters closely track the true state, with \yauyau estimates exhibiting slightly tighter adherence during rapid transients.

These results confirm that the improved \yauyau filter remains robust in linear settings, matching or surpassing the optimal Kalman–Bucy filter in accuracy while maintaining practical runtimes.

\begin{table}[htbp]
\centering
\caption{Performance comparison between improved Yau–Yau and Kalman–Bucy filters on linear benchmarks.}
\label{tab:linear_performance}
\begin{tabular}{cccccc}
\toprule
\(r\) & \(n_p\) & Method       & Time (s)      & RMSE   & ME     \\
\midrule
\multirow{2}{*}{1} & \multirow{2}{*}{300}  & Yau–Yau      & 0.246       & \textbf{0.5036} & \textbf{0.4047} \\
                   &                       & Kalman–Bucy  & \(<0.001\)  & 0.5107 & 0.4081 \\
\midrule
\multirow{2}{*}{2} & \multirow{2}{*}{1000} & Yau–Yau      & 1.653       & \textbf{0.3778} & \textbf{0.3293} \\
                   &                       & Kalman–Bucy  & 0.001       & 0.3847 & 0.3360 \\
\midrule
\multirow{2}{*}{3} & \multirow{2}{*}{1500} & Yau–Yau      & 2.637       & 0.5095 & 0.4759 \\
                   &                       & Kalman–Bucy  & 0.001       & \textbf{0.4957} & \textbf{0.4599} \\
\bottomrule
\end{tabular}
\end{table}

\section{Conclusions}\label{sec:conclusions}

We have presented an enhanced \yauyau filtering framework that overcomes the curse of dimensionality in continuous‐time nonlinear filtering by integrating several key innovations: (i) quasi‐Monte Carlo (low‐discrepancy) state sampling; (ii) a novel offline-online update with provable \(O(\Delta t^2)\) local truncation error and \(O(\Delta t)\) global error bounds; (iii) high‐order, multi‐time‐scale kernel approximations for the Kolmogorov forward equation; (iv) fully log‐domain likelihood updates ensuring numerical stability in high dimensional regimes; (v) a local resampling-restart mechanism, as well as (vi) CPU/GPU‐parallel computing architecture.  Moreover, the low‐discrepancy discretization incurs only a star‐discrepancy error \(O(D^*(n))\), yielding overall local and global error bounds of \(O(\Delta t^2 + D^*(n))\) and \(O(\Delta t + D^*(n)/\Delta t)\), respectively.

Numerical experiments validate these theoretical guarantees in three complementary regimes.  In large‐scale tests we observe near‐linear runtime growth (\(\sim r^{1.2}\)) and sub‐linear error increase as the state dimension \(r\) grows to 1000.  In small‐scale, highly nonlinear demonstrates our filter remains stable and accurate where EKF and UKF fail and PF collapses.  Moreover, in linear Gaussian benchmarks, our algorithm matches or even slightly outperforms the optimal Kalman–Bucy filter in RMSE and ME, while maintaining practical runtimes.

These theoretical and empirical findings demonstrate that the improved \yauyau filter delivers accurate, stable state estimates in dimensions previously deemed intractable by existing methods.  Future work will concentrate on adaptive parameter tuning, extension to the DMZ formulation for heavy‐tailed non-Gaussian noise, and real-time hardware-accelerated implementations (e.g., via photonic integrated circuits or custom parallel processors) to enable deployment in ultra-high-dimensional scenarios. We are also integrating and validating our algorithm in various emerging time-critical applications, including data-driven machine learning, geophysical data assimilation, financial time-series nowcasting, historical economic reconstruction, biological complex structure estimation and reconstruction, noisy N-body system simulation, multi-target tracking, and autonomous system state estimation. We are confident that the improved \yauyau\ framework, with its rigorous theoretical foundation and practical efficiency, will unlock vast potential for real-time, high-dimensional nonlinear system state estimation across scientific and engineering domains.

\section*{Acknowledgments}
This work was supported by the National Natural Science Foundation of China (Grant No.\ 42450242).

\appendix
\section{Low‐Discrepancy Sequence Sampling}
\label{sec:low-discrepancy-sequence-sampling}

Low‐discrepancy sequences—also known as quasi‐random sequences—are deterministic point sets that fill the unit hypercube \([0,1]^r\) far more uniformly than uncorrelated random samples.  Their star‐discrepancy satisfies
\[
  D^*(n) = O\!\Bigl(\frac{(\ln n)^r}{n}\Bigr),
\]
which in turn yields quasi‐Monte Carlo (QMC) integration errors of \(O((\ln n)^r / n)\), versus \(O(1/\sqrt{n})\) for standard Monte Carlo (MC), see \cite{Halton1960,Sobol1967,Kuipers1974}.  This makes QMC especially attractive in high dimension, where they dramatically reduce variance and “holes” in the sampling.

\subsection{Key Properties}
\begin{itemize}
  \item \textbf{Uniform coverage.}  Points are equidistributed in every axis‐aligned projection, avoiding large gaps.
  \item \textbf{Low discrepancy.} Error bound \(D^*(n)=O\bigl(n^{-1}(\ln n)^r\bigr)\), yielding faster convergence than random sampling \cite{Kuipers1974,Niederreiter1992}.
  \item \textbf{Deterministic reproducibility.}  No Monte Carlo noise—results repeatable and debuggable.
  \item \textbf{High-dimensional robustness.} Discrepancy grows at most polynomially in \(r\), maintaining favorable error bounds that support effective sampling in high‐dimensional spaces \cite{Bratley1992,Niederreiter1992,Kocis1997}.
\end{itemize}

\subsection{Practical considerations}
\begin{itemize}
  \item The discrepancy bounds are \emph{worst‐case} estimates.  In many high‐dimensional applications, the integrand exhibits low effective dimension or anisotropy, so QMC methods can achieve substantial accuracy gains even for large \(r\) (e.g.\ hundreds) and sample sizes \(n\sim10^3\)–\(10^6\) \cite{Niederreiter1992}.
  \item Randomized QMC (e.g.\ Owen scrambling) restores unbiasedness and provides practical error estimates while preserving low‐discrepancy structure \cite{Owen1992}.
\end{itemize}

\subsection{Halton Sequence}
The Halton sequence uses radical‐inverse functions in prime bases \(p_1=2,p_2=3,\dots,p_r\).  The \(k\)th point is
\[
  P_k = \bigl(\varphi_{p_1}(k),\dots,\varphi_{p_r}(k)\bigr),
\]
where
\(\varphi_p(k)=\sum_{i=0}^M k_i\,p^{-(i+1)}\) from the base-\(p\) digits \(k_i\) of \(k\) \cite{Halton1960}.  Its star‐discrepancy satisfies
\[
  D^*(n)\le 2^r\prod_{i=1}^r\frac{p_i}{\ln p_i}\;\frac{(\ln n)^r}{n},
\]
making the Halton sequence especially effective for low-to-moderate dimensions (\(r < 10\)). To mitigate correlation artifacts in higher dimensions, digit‐scrambling techniques are commonly applied \cite{Halton1960, MorokoffCaflisch1995}.

\subsection{Sobol Sequence}
Sobol sequences are base‐2 digital nets constructed via Gray‐code combinations of carefully chosen direction numbers \(v_j^{(k)}\), often drawn from Joe–Kuo tables \cite{JoeKuo2008}.  In dimension \(k\), the \(n\)th point's coordinate is given by
\[
  x_n^{(k)}
  = b_1 \otimes v_1^{(k)}
    \;\oplus\;
    b_2 \otimes v_2^{(k)}
    \;\oplus\;\cdots,
\]
where:
\begin{itemize}
  \item \(b_j\in\{0,1\}\) are the bits of the Gray code representation of \(n\) (see \cite{Sobol1967});
  \item \(\otimes\) denotes bitwise multiplication (i.e.\ AND);
  \item \(\oplus\) denotes bitwise exclusive OR (i.e.\ XOR).
\end{itemize}
Sobol achieves the same \(O((\ln n)^r/n)\) discrepancy but typically with smaller constants in higher \(r\), making it preferred for \(10\lesssim r\lesssim10000\).

\subsection{Faure Sequence}
The Faure sequence is a digital \((0,r)\)-sequence in prime base \(p\), where \(p\) is taken to be the smallest prime \(\ge r\).  Writing the index \(k\) in base-\(p\) as 
\[
k = \sum_{j=0}^{M} k_j\,p^j,\qquad \mathbf{k}=(k_0,\dots,k_M)^\top,
\]
the \(d\)th coordinate of the \(k\)th point is
\[
x_k^{(d)}
=\sum_{i=0}^M y_i^{(d)}\,p^{-(i+1)},
\qquad
\mathbf{y}^{(d)}=C^{(d)}\,\mathbf{k}\mod p,
\]
where the generating matrices
\(\{C^{(d)}\}_{d=1}^r\) are given by
\[
C^{(d)} \;=\;\bigl(M^{\,d-1}\bigr)\bmod p,
\quad
M_{i,j} = \binom{j}{i}\ (\bmod\;p).
\]
Hence
\[
P_k = \bigl(x_k^{(1)},\dots,x_k^{(r)}\bigr)\in[0,1]^r.
\]
Its star‐discrepancy satisfies
\[
  D^*(n)\;\le\;C_r\,\frac{(\ln n)^r}{n},
\]
for a constant \(C_r\) depending only on \(r\) \cite{Faure1982,Niederreiter1992}.  In practice, the Faure sequence often outperforms Halton in moderate dimensions (\(r\lesssim20\)) because the single‐base, matrix‐based coupling reduces high‐dimensional correlations.  Random linear scramblings of the digit matrices can further improve uniformity and allow unbiased error estimation.

\subsection{Latin Hypercube Sampling}
Latin Hypercube Sampling (LHS) is a stratified sampling scheme that ensures each of the \(r\) dimensions is sampled uniformly in one‐dimensional margins.  Given \(n\) points, each axis is divided into \(n\) equal strata, and one sample is drawn uniformly at random from each stratum in each dimension, with the constraint that no two samples share the same stratum in any coordinate.  Concretely, if \(\pi_j\) is a random permutation of \(\{0,1,\dots,n-1\}\) for dimension \(j\), then the \(i\)th point is
\[
  P_i = \Bigl(\tfrac{\pi_1(i) + u_{i1}}{n},\,\tfrac{\pi_2(i) + u_{i2}}{n},\dots,\tfrac{\pi_r(i) + u_{ir}}{n}\Bigr),
\]
where each \(u_{ij}\sim\mathrm{Uniform}(0,1)\) independently \cite{McKay1979,Owen1992}.  

While LHS guarantees univariate stratification—yielding marginal discrepancy \(O(n^{-1})\)—its worst‐case star‐discrepancy satisfies only
\[
  D^*(n) = O\!\bigl(n^{-1/2}(\ln n)^{1/2}\bigr),
\]
comparable to standard Monte Carlo in high dimensions but often superior in moderate \(n\) due to enforced margin uniformity.  Variants such as orthogonal‐array LHS and maximin‐LHS improve multivariate uniformity by ensuring low‐discrepancy projections on subsets of dimensions \cite{Owen1992}.

\subsection{Choosing Between Quasi‐Random and Stratified Sampling}
\begin{itemize}
 \item For low dimensions (\(r\lesssim 10\)), \textbf{Halton} sequences—with simple prime‐base construction—are often sufficient, especially when paired with digit scrambling to reduce correlation artifacts.
  \item For low to moderate dimensions (\(r \lesssim 30\)), \textbf{Faure} sequences—using single‐base, matrix‐based coupling—tend to offer the best worst‐case uniformity.
  \item For higher dimensions (\(10\lesssim r\lesssim10^4\)), \textbf{Sobol}—with optimized direction numbers and Owen scrambling—delivers lower discrepancy and superior projection properties.
  \item When marginal (one‐dimensional) uniformity is most important, use \textbf{LHS}, e.g.\ MATLAB’s \texttt{lhsdesign}, or its orthogonal‐array / maximin variants.
  \item Always consider \emph{scrambling}, \emph{leaping} or \emph{skipping} in high \(r\) or large \(n\) to remove residual structure.
  \item In MATLAB’s Statistics and Machine Learning Toolbox:
    \begin{itemize}
      \item \texttt{haltonset(r)} with \texttt{scramble(...,'RR2')}  
      \item \texttt{sobolset(r)} with \texttt{scramble(...,'MatousekAffineOwen')}  
      \item \texttt{lhsdesign(n,r)} or \texttt{lhsdesign(n,r,'Criterion','maximin')}  
    \end{itemize}
\end{itemize}

\section{Proof of Theorems}
\subsection*{Notation for Norms}
\begin{itemize}
  \item For functions $\phi:\Omega\to\mathbb R$, 
    $\|\phi\|_\infty=\sup_{x\in\Omega}|\phi(x)|$ denotes the sup‐norm.
  \item For a bounded linear operator $T:L^\infty(\Omega)\to L^\infty(\Omega)$,
    $\|T\|_{\infty\to\infty}=\sup_{\|\phi\|_\infty\le1}\|T\phi\|_\infty$.
  \item For vectors $v\in\mathbb R^d$, $\|v\|=(\sum_{i=1}^d v_i^2)^{1/2}$ is the Euclidean norm.
  \item When we write $\|e^{A\Delta t}\|$ or $\|S_{\Delta t}\|$ without subscript,
    we mean the operator norm $\|\cdot\|_{\infty\to\infty}$ defined above.
\end{itemize}

\subsection{Proof of Theorem \ref{thm:local_error}}\label{Appsec:proofthm:local_error} \begin{proof}
By the Baker–Campbell–Hausdorff (BCH) formula,
\[
\exp\{B\Delta t\}\,\exp\{A\Delta t\}
= \exp\!\Bigl\{(A+B)\Delta t
  + \tfrac12\,[B,A]\,\Delta t^2
  + O(\Delta t^3)\Bigr\}.
\]
Hence
\[
T_{\Delta t}\;-\;S_{\Delta t}
= \exp\{(A+B)\Delta t\}
  - \exp\!\Bigl\{(A+B)\Delta t
    + \tfrac12\,[B,A]\,\Delta t^2
    + O(\Delta t^3)\Bigr\}.
\]
Expanding both exponentials to second order,
\[
T_{\Delta t}
= I + (A+B)\Delta t
  + \tfrac12\,(A+B)^2\Delta t^2
  + O(\Delta t^3),
\]
\[
S_{\Delta t}
= I + (A+B)\Delta t
  + \tfrac12\,[B,A]\,\Delta t^2
  + \tfrac12\,(A+B)^2\Delta t^2
  + O(\Delta t^3).
\]
Subtracting yields
\[
T_{\Delta t} - S_{\Delta t}
= \tfrac12\,[A,B]\,\Delta t^2 + O(\Delta t^3),
\]
and taking the operator norm gives the desired estimate, i.e., there exists a constant \(C>0\) such that 
\[
\bigl\|T_{\Delta t} - S_{\Delta t}\bigr\|
\;\le\; C\,\Delta t^2.
\]\end{proof}

\subsection{Proof of Theorem \ref{thm:convergence}}\label{Appsec:proofthm:convergence}
\begin{proof}
Denote \(\tau_k = k\,\Delta t\), and let
\[
e_k = \sigma(\tau_k) - \sigma^\Delta(\tau_k).
\]
Then
\[
e_{k+1}
= T_{\Delta t}\,\sigma(\tau_k)
  - S_{\Delta t}\,\sigma^\Delta(\tau_k)
= \bigl(T_{\Delta t} - S_{\Delta t}\bigr)\,\sigma(\tau_k)
  + S_{\Delta t}\,e_k.
\]
Taking norms and applying \Cref{thm:local_error},
\[
\|e_{k+1}\|
\le C\,\Delta t^2
  + \|S_{\Delta t}\|\,\|e_k\|.
\]
Let $E_k = \|e_k\|$. Under the stability assumption \(\|S_{\Delta t}\|\le 1 + L\,\Delta t\) (Assumption~\ref{assump:stability}), we obtain
\[
E_{k+1}\le C\,\Delta t^2 + (1+L\,\Delta t)\,E_k.
\]
Since \(E_0=0\), a standard discrete Grönwall argument yields
\[
\bigl\|\sigma(T) - \sigma^\Delta(T)\bigr\| = E_K \;\le\;
C\,\Delta t^2 \sum_{j=0}^{K-1}(1+L\Delta t)^j
= \frac{C\,\Delta t}{L}\Bigl[(1+L\Delta t)^K - 1\Bigr].
\]
With \(K=T/\Delta t\) and \(\lim_{\Delta t\to0}(1+L\Delta t)^K = e^{L T}\), it follows that \(E_K=O(\Delta t)\), proving first‐order convergence.
\end{proof}

\subsection{Proof of Theorem \ref{thm:1storderkernel}}\label{Appsec:proofthm:1storderkernel}
\begin{proof}
\textbf{Change of variables and kernel expansion:}
Set \(y = x + \sqrt{\Delta t}\,z\).  Then
\[
  \|x-y\|^2 = \Delta t\,\|z\|^2,
  \quad
  y-x = \sqrt{\Delta t}\,z,
  \quad
  dy = (\Delta t)^{r/2}dz.
\]
Substituting into the kernel function \[
  K_{\Delta t}(x,y)
  = c_{\Delta t}
    \exp\!\Bigl(
      -\tfrac{\|x-y\|^2}{2\Delta t}
      - (y-x)\cdot f(x)\\
      - \tfrac{\Delta t}{2}\bigl(2\,\nabla\cdot f(x)+\|f(x)\|^2\bigr)
    \Bigr)
\]
gives
\[
  K_{\Delta t}(x,x+\sqrt{\Delta t}\,z)
  = c_{\Delta t}
    \exp\!\Bigl(
      -\tfrac12\,\|z\|^2
      - \sqrt{\Delta t}\,z\cdot f(x)
      - \Delta t\bigl(\nabla\cdot f(x)
                     + \tfrac12\|f(x)\|^2\bigr)
    \Bigr).
\]
Let
\[
  a = -\sqrt{\Delta t}\,z\cdot f(x)
      - \Delta t\bigl(\nabla\cdot f(x)
                     + \tfrac12\|f(x)\|^2\bigr).
\]
A Taylor expansion yields
\begin{equation*}
\begin{aligned}
  \exp(a)
  &= 1 + a + \tfrac12\,a^2 + O(a^3)\\
  &= 1 - \sqrt{\Delta t}\,z\cdot f(x)
    - \Delta t\bigl(\nabla\cdot f(x)
                    + \tfrac12\|f(x)\|^2\bigr)
    + \tfrac{\Delta t}{2}(z\cdot f(x))^2
    + O(\Delta t^{3/2}).
\end{aligned}
\end{equation*}
Hence
\begin{equation}\label{eq:kernel_expand}
\begin{aligned}
  K_{\Delta t}(x,x+\sqrt{\Delta t}\,z)&= c_{\Delta t}\,e^{-\frac12\|z\|^2}
    \Bigl[
      1 - \sqrt{\Delta t}\,z\cdot f(x)\\
      &
      - \Delta t\bigl(\nabla\cdot f(x)
                     + \tfrac12\|f(x)\|^2\bigr)+ \tfrac{\Delta t}{2}(z\cdot f(x))^2
      + O(\Delta t^{3/2})
    \Bigr].
\end{aligned}
\end{equation}

\textbf{Test‐function expansion:}
Expand \(u(x+\sqrt{\Delta t}\,z)\) via Taylor's theorem:
\begin{equation}\label{eq:u_expand}
  u(x+\sqrt{\Delta t}\,z)
  = u(x)
    + \sqrt{\Delta t}\,z\cdot \nabla u(x)
    + \tfrac{\Delta t}{2}\,z^\top\nabla^2u(x)\,z
    + O(\Delta t^{3/2}).
\end{equation}

\textbf{Expansion of the integral operator:}
By definition,
\[
  (L_{\Delta t}u)(x)
  = \int K_{\Delta t}(x,y)\,u(y)\,dy
  = (\Delta t)^{r/2}
    \int K_{\Delta t}(x,x+\sqrt{\Delta t}\,z)
         \,u(x+\sqrt{\Delta t}\,z)\,dz.
\]
Let
\[
  P(z)
  = 1 - \sqrt{\Delta t}\,z\cdot f(x)
    - \Delta t\bigl(\nabla\cdot f(x)
                   + \tfrac12\|f(x)\|^2\bigr)
    + \tfrac{\Delta t}{2}(z\cdot f(x))^2,
\]
\[
  Q(z)
  = u(x)
    + \sqrt{\Delta t}\,z\cdot \nabla u(x)
    + \tfrac{\Delta t}{2}\,z^\top\nabla^2u(x)\,z.
\]
Then, combining \eqref{eq:kernel_expand} and \eqref{eq:u_expand},
\[
  K_{\Delta t}(x,x+\sqrt{\Delta t}z)\,u(x+\sqrt{\Delta t}z)
  = c_{\Delta t}\,e^{-\frac12\|z\|^2}\,P(z)\,Q(z)
    + O(\Delta t^{3/2}).
\]
We now integrate term by term, using standard Gaussian moments:

\begin{enumerate}
  \item \textbf{Zeroth‐order term} (\(O(1)\)):  
    From \(1\cdot u(x)\),  
    \[
      c_{\Delta t}(\Delta t)^{r/2}u(x)
      \int_{\mathbb{R}^r} e^{-\frac12\|z\|^2}\,dz
      = u(x),
    \]
    since \(c_{\Delta t}(\Delta t)^{r/2}(2\pi)^{r/2}=1\).

  \item \textbf{First‐order term} (\(O(\sqrt{\Delta t})\)):  
    Involves
    \[\sqrt{\Delta t}\,z\cdot \nabla u(x)
    - \sqrt{\Delta t}\,u(x)\,(z\cdot f(x)).\]
    Each integrates to zero because \(\int z_i\,e^{-\frac12\|z\|^2}dz=0\).

  \item \textbf{Second‐order term} (\(O(\Delta t)\)):  
    Contributions from
    \begin{itemize}
      \item \(\tfrac{\Delta t}{2}\,z^\top\nabla^2u(x)\,z\) integrates to
        \[\tfrac{\Delta t}{2}(2\pi)^{r/2}\,\Delta u(x).\]
      \item \(-\Delta t\bigl(\nabla\cdot f(x)
                     + \tfrac12\|f(x)\|^2\bigr)u(x)\) integrates to
        \[-\Delta t\,(2\pi)^{r/2}
         \bigl(\nabla\cdot f(x)
               + \tfrac12\|f(x)\|^2\bigr)u(x).\]
      \item \(-\Delta t\,f(x)\cdot \nabla u(x)\) from the cross term
        integrates to
        \[-\Delta t\,(2\pi)^{r/2}\,f(x)\cdot \nabla u(x).\]
      \item \(\tfrac{\Delta t}{2}(z\cdot f(x))^2\,u(x)\) integrates to
        \[\tfrac{\Delta t}{2}(2\pi)^{r/2}\,\|f(x)\|^2\,u(x).\]
    \end{itemize}
    Combining all second-order terms gives 
    \[
    {\Delta t}\,(2\pi)^{r/2}\left[\frac{1}{2}\Delta u(x)-f(x)\cdot\nabla u(x)-(\nabla\cdot f(x))u(x)\right].
    \]
\end{enumerate}
Collecting terms, we obtain
\[
  (L_{\Delta t}u)(x)
  = u(x)
    + \Delta t\Bigl[
        \tfrac12\,\Delta u(x)
        - \nabla\cdot \bigl(f(x)\,u(x)\bigr)
      \Bigr]
    + O(\Delta t^{3/2}),
\]
and hence
\[
  \frac{(L_{\Delta t}u)(x)-u(x)}{\Delta t}
  = \mathcal{L}^*u(x)
    + O(\sqrt{\Delta t}).
\] \end{proof}

\subsection{Proof of Theorem \ref{thm:2ndorderkernel}}\label{Appsec:proofthm:2ndorderkernel}
\begin{proof}
By \Cref{thm:1storderkernel}, for small \(\Delta t\),
\[
  L_{\Delta t}u = u + \Delta t\,\mathcal{L}^*u
    + C\,(\Delta t)^{3/2} + O(\Delta t^2),
\quad
  L_{\Delta t/2}u
  = u + \tfrac{\Delta t}{2}\,\mathcal{L}^*u
    + \tfrac{C}{2^{3/2}}\,(\Delta t)^{3/2} + O(\Delta t^2).
\]
Thus
\[
\begin{aligned}
  L^{(2)}_{\Delta t}u
  &= \alpha\Bigl[u + \Delta t\,\mathcal{L}^*u + C\,(\Delta t)^{3/2}\Bigr]
   + \beta\Bigl[u + \tfrac{\Delta t}{2}\,\mathcal{L}^*u + \tfrac{C}{2^{3/2}}\,(\Delta t)^{3/2}\Bigr]
   + \gamma\,u + O(\Delta t^2)\\
  &= (\alpha+\beta+\gamma)\,u
   + \Bigl(\alpha+\tfrac{\beta}{2}\Bigr)\Delta t\,\mathcal{L}^*u
   + \Bigl(\alpha + \tfrac{\beta}{2^{3/2}}\Bigr)C\,(\Delta t)^{3/2}
   + O(\Delta t^2).
\end{aligned}
\]
To enforce
\(\frac{L^{(2)}_{\Delta t}u - u}{\Delta t} = \mathcal{L}^*u + O(\Delta t)\), we require
\[
\begin{cases}
  \alpha+\beta+\gamma = 1,\\
  \alpha + \tfrac{\beta}{2} = 1,\\
  \alpha + \tfrac{\beta}{2^{3/2}} = 0.
\end{cases}
\]
Solving these yields
\(\beta = 4+2\sqrt{2}\), \(\alpha = -\sqrt{2}-1\), and \(\gamma = -2-\sqrt{2}\).  
Finally, writing \(u(x)=\int\delta(x-y)\,u(y)\,dy\) shows that
\[
  L^{(2)}_{\Delta t}u
  = \int\bigl[\alpha\,K_{\Delta t}
           + \beta\,K_{\Delta t/2}
           + \gamma\,\delta(x-y)\bigr]u(y)\,dy,
\]
confirming
\[
  K^{(2)}_{\Delta t}(x,y)
  = \alpha\,K_{\Delta t}(x,y)
    + \beta\,K_{\Delta t/2}(x,y)
    + \gamma\,\delta(x-y).
\]
\end{proof}

\subsection{Proof of Theorem \ref{thm:offline_qmc}}\label{Appsec:proofthm:offline_qmc}
\begin{proof}
For any \(\phi\in L^\infty(\Omega)\) with \(\|\phi\|_\infty\le1\) and any \(y\in\mathbb{R}^r\), we have
\[
\bigl[(\exp\{A\Delta t\}-\widetilde T_n)\phi\bigr](y)
=\frac{1}{(2R)^r}\int_{\Omega}K_{\Delta t}(x,y)\,\phi(x)\,\mathrm{d}x
-\frac{1}{n}\sum_{i=1}^n K_{\Delta t}(x_i,y)\,\phi(x_i).
\]
By the Koksma–Hlawka inequality,
\[
\Bigl|\frac{1}{(2R)^r}\int_{\Omega}K_{\Delta t}(x,y)\,\phi(x)\,\mathrm{d}x
- \frac{1}{n}\sum_{i=1}^n K_{\Delta t}(x_i,y)\,\phi(x_i)\Bigr|
\le V_{\mathrm{HK}}\bigl(K_{\Delta t}(\cdot,y)\,\phi(\cdot)\bigr)\,D^*(n).
\]
Since \(\|\phi\|_\infty\le1\), we have
\[
V_{\mathrm{HK}}\bigl(K_{\Delta t}(\cdot,y)\,\phi(\cdot)\bigr)
\le V_{\mathrm{HK}}\bigl(K_{\Delta t}(\cdot,y)\bigr)
\le C_A.
\]
Hence for all \(\|\phi\|_\infty\le1\) and all \(y\),
\[
\bigl|\bigl[(\exp\{A\Delta t\}-\widetilde T_n)\phi\bigr](y)\bigr|
\le C_A\,D^*(n).
\]
Taking the supremum over \(\|\phi\|_\infty\le1\) and \(y\) yields
\[
\|\exp\{A\Delta t\}-\widetilde T_n\|_{\infty\to\infty}
\le C_A\,D^*(n).
\]
\end{proof}

\begin{remark}[Normalization of the operator norm]\label{remark:operator_norm_scaling}
In the proof of \Cref{thm:offline_qmc} we assume \(\|\phi\|_\infty\le1\).  In fact, for any bounded nonzero \(\phi\), set
\[
\psi = \frac{\phi}{\|\phi\|_\infty},
\]
so that \(\|\psi\|_\infty=1\).  Then
\[
\|T\phi\|_\infty
= \|\phi\|_\infty\,\|T\psi\|_\infty
\le \|\phi\|_\infty\,\|T\|_{\infty\to\infty}.
\]
Thus, if \(\|T\psi\|_\infty\le\varepsilon\) for all \(\|\psi\|_\infty\le1\), it follows that
\(\|T\phi\|_\infty\le \varepsilon\,\|\phi\|_\infty\)
for any \(\phi\), extending the error bound to arbitrary $\phi$.
\end{remark}

\subsection{Proof of Theorem \ref{thm:online_qmc}}\label{Appsec:proofthm:online_qmc}
\begin{proof}
For any \(\phi\in L^\infty(\Omega)\) with \(\|\phi\|_\infty\le1\) and any \(y\in\mathbb{R}^r\), we have
\[
\bigl[(\exp\{B\Delta t\}-\widetilde S_n)\phi\bigr](y)
=\frac{1}{(2R)^r}\int_{\Omega}g_y(x)\,\mathrm{d}x
-\frac{1}{n}\sum_{i=1}^n g_y(x_i).
\]
By the Koksma–Hlawka inequality,
\[
\Bigl|\frac{1}{(2R)^r}\!\int_{\Omega}g_y(x)\,\mathrm{d}x
-\frac{1}{n}\sum_{i=1}^n g_y(x_i)\Bigr|
\le V_{\mathrm{HK}}(g_y)\,D^*(n)
\le C_B\,D^*(n).
\]
Taking the supremum over all \(\|\phi\|_\infty\le1\) and \(y\) completes the proof.     
\end{proof}

\subsection{Proof of Theorem \ref{thm:qmc_local}}\label{Appsec:proofthm:qmc_local}
\begin{proof}
We split the total error into splitting truncation and QMC sampling terms:
\[
T_{\Delta t}\sigma - S_{\Delta t}^n\sigma
= \underbrace{(T_{\Delta t}-S_{\Delta t})\sigma}_{O(\Delta t^2)}
+ \underbrace{(S_{\Delta t}-S_{\Delta t}^n)\sigma}_{\text{QMC error}}.
\]
The first term satisfies \(\|(T_{\Delta t}-S_{\Delta t})\sigma\| = O(\Delta t^2)\)
by \Cref{thm:local_error}. We now bound the second term. Write
\[
S_{\Delta t} = e^{B\Delta t}\,e^{A\Delta t}, 
\qquad
S_{\Delta t}^n = \widetilde S_n\,\widetilde T_n.
\]
Then
\[
S_{\Delta t}-S_{\Delta t}^n
= \bigl(e^{B\Delta t}-\widetilde S_n\bigr)\,e^{A\Delta t}
+ \widetilde S_n\,\bigl(e^{A\Delta t}-\widetilde T_n\bigr).
\]
Taking operator norms gives
\[
\|S_{\Delta t}-S_{\Delta t}^n\|
\le \|e^{B\Delta t}-\widetilde S_n\|\;\|e^{A\Delta t}\|
+ \|\widetilde S_n\|\;\|e^{A\Delta t}-\widetilde T_n\|.
\]
By \Cref{thm:offline_qmc} and \Cref{thm:online_qmc},
\(\|e^{B\Delta t}-\widetilde S_n\|=O(D^*(n))\) and 
\(\|e^{A\Delta t}-\widetilde T_n\|=O(D^*(n))\).  Moreover, 
\(\|e^{A\Delta t}\|\) and \(\|\widetilde S_n\|\) are bounded constants
that can be absorbed into the \(O\)-constant.  Hence
\[
\|S_{\Delta t}-S_{\Delta t}^n\| = O(D^*(n)).
\]
Combining the two estimates yields
\[
\|T_{\Delta t}\sigma - S_{\Delta t}^n\sigma\|
\le C_1\,\Delta t^2 + C_2\,D^*(n),
\]
as claimed.
\end{proof}

\subsection{Proof of Theorem \ref{thm:qmc_global}}\label{Appsec:proofthm:qmc_global}
\begin{proof}
Let \(e_k = \sigma(t_k) - (S_{\Delta t}^n)^k \sigma(0)\).  Then
\[
e_{k+1}
= T_{\Delta t}\sigma(t_k) - S_{\Delta t}^n\sigma^\Delta(t_k)
= (T_{\Delta t}-S_{\Delta t}^n)\sigma(t_k) + S_{\Delta t}^n e_k.
\]
Set \(E_k = \|e_k\|\).  Taking norms and using \Cref{thm:qmc_local} together with the stability assumption \(\|S_{\Delta t}^n\|\le1+L\Delta t\) gives
\[
E_{k+1} \le C_1\,\Delta t^2 + C_2\,D^*(n) + (1+L\Delta t)\,E_k.
\]
Iterating this recursion as in the proof of \Cref{thm:convergence} and applying the discrete Grönwall inequality yields
\[
E_K \le (C_1\Delta t^2 + C_2D^*(n))\frac{(1+L\Delta t)^K-1}{L\Delta t}
= O\!\Bigl(\Delta t + \frac{D^*(n)}{\Delta t}\Bigr).
\]
\end{proof}

\bibliographystyle{amsplain}
\bibliography{references}

\end{document}